\documentclass{article}[11pt]
\usepackage{curves,epic,epsfig,makeidx,amsmath,amsfonts,latexsym,subfigure}
\newtheorem{theorem}{Theorem}[section]
\newtheorem{prop}{Proposition}[section]
\newtheorem{lemma}{Lemma}[section]
\newtheorem{conj}{Conjecture}[section]

\begin{document}

\title{\large\textbf{{A STOCHASTIC MODEL FOR COMPETING GROWTH ON $\mathbb{R}^d$}}}

\author{Maria Deijfen \thanks{Stockholm University. E-mail: mia@matematik.su.se}
\\ Olle H\"{a}ggstr\"{o}m \thanks{Chalmers University
of Technology. E-mail: olleh@math.chalmers.se}\\
Jonathan Bagley \thanks{UMIST, Manchester. E-mail:
jonathan.bagley@umist.ac.uk}}

\date{2003}

\maketitle

\thispagestyle{empty}

\begin{abstract}

\noindent A stochastic model, describing the growth of two
competing infections on $\mathbb{R}^d$, is introduced. The growth
is driven by outbursts in the infected region, an outburst in the
type 1 (2) infected region transmitting the type 1 (2) infection
to the previously uninfected parts of a ball with stochastic
radius around the outburst point. The main result is that with the
growth rate for one of the infection types fixed, mutual unbounded
growth has probability zero for all but at most countably many
values of the other infection rate. This is a continuum analog of
a result of H\"{a}ggstr\"{o}m and Pemantle. We also extend a shape
theorem of Deijfen for the corresponding model with just one type
of infection.

\vspace{1cm}


\vspace{0.5cm}

\hspace*{4.8cm} Secondary 82B43.
\end{abstract}

\section{Introduction}

\noindent In Deijfen (2003), a model is introduced that describes
the random growth of an infected region in $\mathbb{R}^d$ by means
of spherical outbursts in the infected region. The purpose of the
present paper is to extend the work of Deijfen (2003) in two
directions. First, we generalize the asymptotic shape theorem of
that paper from bounded outburst radii to unbounded ones
satisfying a certain moment condition; see Theorem
\ref{th:aformgen}. Second, we extend the model to encompass two
competing types of infection, and prove a continuum analog of a
result of H\"{a}ggstr\"{o}m and Pemantle (2000) concerning the
impossibility of mutual unbounded growth; see Theorem
\ref{th:huvudres}.\medskip

\noindent The model in Deijfen (2003) can be viewed as a
generalization to continuous space of the well-known Richardson
model. The Richardson model, first introduced in Richardson
(1973), describes growth on $\mathbb{Z}^d$. Sites can be either
infected or uninfected: An uninfected site becomes infected at a
rate proportional to the number of infected nearest neighbors and
once infected it never recovers. In the continuum model, the
growth takes place by way of outbursts in the infected region, an
outburst at an infected point causing a ball with stochastic
radius around the outburst point to be infected. Hence, for all
$t$, the infected region at time $t$, denoted by $S_t$, is a union
of randomly sized Euclidean balls. The dynamics is that, given
$S_t$, the time until the next outburst occurs is exponentially
distributed with expected value proportional to $|S_t|^{-1}$ and
the outburst point is chosen uniformly in $S_t$. The main result
in Deijfen (2003) is a shape theorem which asserts that, if there
is an upper bound for the radii of the outburst balls, then on the
scale $1/t$ the set $S_t$ has an asymptotic shape, which due to
rotational invariance must be a Euclidean ball. The following
result states that the conclusion of the shape theorem in Deijfen
(2003) is still valid under the weaker assumption that the radius
distribution, denoted by $F$, has a moment generating function.
Here $B(x,r)$ is a closed ball with radius $r$ centered at $x$.

\begin{theorem}[Generalized shape theorem]\label{th:aformgen} Fix
$d\geq 1$ and consider the $d$-dimensional continuum growth model
with rate $\lambda$. Assume that

\begin{equation}\label{eq:rvillk}
\int_0^\infty e^{-\varphi r}dF(r)<\infty \textrm{ for some }
\varphi<0
\end{equation}

\noindent and let $S_0\subset \mathbb{R}^d$ be arbitrary but
bounded with strictly positive Lebesgue measure. Then there is a
real number $\mu>0$ such that, for any
$\varepsilon\in(0,\lambda{\mu}^{-1})$, almost surely
$$
(1-\varepsilon)B\left(0,\lambda{\mu}^{-1}\right)\subset
\frac{S_t}{t}
\subset(1+\varepsilon)B\left(0,\lambda{\mu}^{-1}\right)
$$
for all sufficiently large $t$.
\end{theorem}

\noindent \textbf{Remark 1.1} The condition (\ref{eq:rvillk}) is
probably not optimal. It has been suggested by a referee that it
should be possible to replace with a $p$-moment assumption for
$p>2d$.\medskip

\noindent In H\"aggstr\"om and Pemantle (1998) a two-type version
of the Richardson model is introduced: Two particle types -- type
1 and type 2 -- compete for space on $\mathbb{Z}^d$, the dynamics
being that an empty site becomes occupied by a type $i$ particle
at a rate that is proportional to the number of nearest type $i$
neighbors and an occupied site remains occupied forever. In this
paper a similar two-type version of the continuum model is
introduced. A careful description will be given in Section 3 but
roughly the model is as follows: Two different non-reversible
entities -- which henceforth will be referred to as type 1 and
type 2 infection -- compete for space on $\mathbb{R}^d$. The
infected region at time $t$ can be divided in two disjoint sets
$S_t^1$ and $S_t^2$ indicating the region occupied by type 1
infection and type 2 infection respectively. As in the one-type
model the growth takes place by way of outbursts which infect the
previously uninfected parts of a ball with random radius around
the outburst point. The outbursts are of two types: Type 1
outbursts occur in the type 1 infected region and result in
outburst balls of type 1, that is, the infection transmitted by
the outburst is of type 1. Type 2 outbursts occur in the type 2
infected region and causes outburst balls loaded with type 2
infection. The radii of the outburst balls are i.i.d.\ random
variables with the same distribution, denoted by $F$, for both
outburst types. For $i=1,2$, given $S_t^i$ the time until an
outburst occurs in $S_t^i$ is exponentially distributed with
parameter $\lambda_i|S_t^i|$ and the outburst point is uniformly
distributed over $S_t^i$. Note that, if $\lambda_1=\lambda_2$ and
if we do not distinguish between the infection types, this model
is equivalent to the one-type model introduced in Deijfen
(2003).\medskip

\noindent Consider the development of the infection in the
two-type Richardson model. There are two possible scenarios:

\begin{itemize}
\item[1.] One of the infection types is at some point surrounded by the
other type, implying that only finitely many sites are ever
occupied by the surrounded type.
\item[2.] Both infection types keep growing indefinitely.
\end{itemize}

\noindent Clearly the first scenario has positive probability, but
what about the second one? This issue is dealt with in
H\"aggstr\"om and Pemantle (1998) and (2000). The main result in
the first paper is that, if $\lambda_1=\lambda_2$ -- that is, if
the infection types are equal in power -- and $d=2$, then the
event $G=\{$both infection types grow indefinitely$\}$ has
positive probability. In the second paper the case $\lambda_1\neq
\lambda_2$ is considered and the main result is that, if
$\lambda_1=1$, then $G$ has probability zero for all but at most
countably many values of $\lambda_2$ (note that, by time-scaling,
the assumption that $\lambda_1=1$ is no restriction). The main
result in the present paper is that, for $d\geq 2$, this holds
also in the two-type continuum model. This strongly suggests that
$G$ in fact has probability zero for \emph{all} choices of
$\lambda_2\neq 1$. See H\"aggstr\"om and Pemantle (2000) for some
intuitive reasoning behind this statement.\medskip

\noindent Before formulating the result, we have to specify what
the event ``both infection types grow indefinitely'' will mean in
the continuum model. To this end, let
$$
G_i=\{\textrm{the type $i$ infection reaches points arbitrarily
far from the origin}\}
$$
and define $G=G_1\cap G_2$ so that $G$ hence is the event that
both infection types reach arbitrarily far away from the origin
simultaneously. If the outburst radius is bounded, the type 1 (2)
infection is clearly prevented from growing any further if it is
surrounded by a type 2 (1) layer whose thickness exceeds the upper
bound for the size of an outburst. In this case the event $G$ thus
means that none of the infection types is enclosed by the other
and so the concept of co-existence is similar to the lattice case.
On the other hand, when the outburst radius is not bounded, it is
not possible to rule out $G$ in finite time, indicating that we
are dealing with a more subtle concept of co-existence.\medskip

\noindent It is not hard to show that the events $G_1$ and $G_2$
both have positive probability for all choices of $\lambda_1$ and
$\lambda_2$ (see Proposition \ref{prop:pgisn}). The event
$G=G_1\cap G_2$ is more complicated to study. However, intuitively
it is clear that $G$ should not occur if $\lambda_1\neq
\lambda_2$: That both infection types reach arbitrarily far away
from the origin means that some kind of balance of power reigns
between the infection types and if one of them is more powerful
than the other there is no reason to believe that such a balance
should be possible. Our main result is a step on the way towards a
confirmation of this intuition. To formulate it some notation is
needed. Let $P_{\Gamma_1,\Gamma_2}^{\lambda_1,\lambda_2}$ denote
the probability law of the two-type growth process started at time
zero from $S_0^1=\Gamma_1$ and $S_0^2=\Gamma_2$ and with infection
rates $\lambda_1$ and $\lambda_2$ respectively. Our first result
for the two-type model is that, under the assumption that
arbitrarily small outbursts are possible, the particular choice of
$\Gamma_1$ and $\Gamma_2$ is irrelevant in deciding whether the
event $G$ has positive probability or not.

\begin{prop}\label{prop:startomr}
Let $(\Gamma_1,\Gamma_2)$ and $(\Gamma_1',\Gamma_2')$ be two pairs
of disjoint, bounded subsets of $\mathbb{R}^d$ with strictly
positive Lebesgue measures. Furthermore, suppose that $F$ has
unbounded support and satisfies $F(\varepsilon)>0$ for all
$\varepsilon>0$. Then
$$
P_{\Gamma_1,\Gamma_2}^{\lambda_1,\lambda_2}(G)>0\Rightarrow
P_{\Gamma_1',\Gamma_2'}^{\lambda_1,\lambda_2}(G)>0.
$$
\end{prop}

\noindent \textbf{Remark 1.2} The proof of Proposition 1.1 can
easily be modified to cover the case with bounded support as well,
provided the following (obviously necessary) condition on
$(\Gamma'_1, \Gamma'_2)$: If the radius distribution is bounded by
$r$, then neither of $\Gamma'_1$ or $\Gamma'_2$ may contain an
impenetrable layer of thickness $r$ around the other.\medskip

\noindent In words, if we can find two sets $\Gamma_1$ and
$\Gamma_2$ such that the event $G$ has positive probability when
starting from $\Gamma_1$ and $\Gamma_2$, then it follows that $G$
has positive probability for all other initial sets $\Gamma_1'$
and $\Gamma_2'$ as well. Thus we may restrict our attention to the
case when the process is started from two balls with radius given
by the mean outburst radius, denoted by $\gamma$. The balls are
placed next to each other, one of them centered at the point
$-\textbf{2}\gamma=(-2\gamma,0,\ldots ,0)$ and the other at the
origin. The notation for the probability law in this case is
simplified by dropping the subscripts so that
$P^{\lambda_1,\lambda_2}$ hence denotes the law of the growth
process started from $\Gamma_1=B(-\textbf{2}\gamma,\gamma)$ and
$\Gamma_2=B(0,\gamma)$. Furthermore, note that by time scaling we
may assume that $\lambda_1=1$. The main result is as follows.

\begin{theorem}\label{th:huvudres}
If $F$ satisfies (\ref{eq:rvillk}), then, for $d\geq 2$, the set
$\{\lambda_2;\hspace{0.1cm} P^{1,\lambda_2}(G)>0\}$ is countable.
\end{theorem}

\noindent Thus $P^{1,\lambda_2}(G)=0$ for all but at most
countably many values of $\lambda_2$. As mentioned before, this
strongly suggests that $P^{1,\lambda_2}(G)=0$ for all
$\lambda_2\neq 1$. In the case $\lambda_2=1$ on the other hand, it
seems reasonable to suspect that $G$ has positive probability.
That $\lambda_2=1$ means that the two infection types are equal in
power and hence it should be possible for some kind of balance of
power to arise. Let us summarize all this in the following
conjecture.\medskip

\begin{conj}
If (\ref{eq:rvillk}) holds, we have $\{\lambda_2;\hspace{0.1cm}
P^{1,\lambda_2}(G)>0\}=\{1\}$ for $d\geq 2$.
\end{conj}

\noindent The rest of the paper is organized as follows. Theorem
\ref{th:aformgen} is proved in Section 2. Section 3 provides a
closer description of the two-type continuum model and Section 4
contains a number of auxiliary results needed in the following
sections. In Section 5 an analogue of the ``key proposition'' in
H\"aggstr\"om and Pemantle (2000) is formulated and proved. This
result will be of vital importance in the proof of Theorem
\ref{th:huvudres} and due to the fact that the asymptotic shape
for the continuum model is known to be a Euclidean ball, its proof
is somewhat more appealing for the intuition compared to the
lattice case. Theorem \ref{th:huvudres} is proved in Section 6 and
Proposition \ref{prop:startomr} in Section 7.

\section{The asymptotic shape}

\noindent The aim in this section is to prove Theorem
\ref{th:aformgen}.\medskip

\noindent As described in the introduction, the growth in the
one-type process is generated by stochastically sized spherical
outbursts in the infected region $S_t$. Given the development of
the infection up to time $t$, the time until an outburst occurs
somewhere in $S_t$ is exponentially distributed with parameter
$\lambda|S_t|$ and the location of the outburst is uniformly
distributed over $S_t$. To formally construct the model in $d$
dimensions, a $(d+1)$-dimensional Poisson process with rate
$\lambda$ is used, the extra dimension representing time. A
bounded set $\Gamma\subset \mathbb{R}^d$ with strictly positive
Lebesgue measure is picked to initiate the growth and also, to the
points in the Poisson process, i.i.d.\ radius variables with
distribution $F$ are attached. Starting at time zero the growth is
then brought about by following the cylinder
$\Gamma\times\mathbb{R}$ upwards along the time axis until a point
in the Poisson process is found. An outburst then takes place at
this point and an infection ball $B_1$ with radius given by the
radius variable associated with this particular Poisson point is
created around the outburst point. The infected region is now
given by $\Gamma\cup B_1$. Scanning within the cylinder
$(\Gamma\cup B_1)\times \mathbb{R}$ further upwards along the time
axis, a new Poisson point is eventually hit and a new infection
ball $B_2$ arises around this point. And so on. For a more
thorough description of the construction we refer to Deijfen
(2003).\medskip

\noindent The main result in Deijfen (2003) is the asymptotic
shape result in Theorem \ref{th:aformgen} under the stronger
assumption that the outbursts radii are bounded. (The result was
formulated for $\lambda=1$ only, but the general result follows by
a simple time scaling argument.)\medskip

\noindent Let $\tilde{T}(x)$ be the time when the entire ball with
radius $\gamma$ around the point $x$ is infected in a unit rate
process started from $S_0=B(0,\gamma)$. The time constant $\mu$ is
given by

\begin{equation}\label{eq:mu}
\mu=\lim_{n\rightarrow \infty}
\frac{\textrm{E}[\tilde{T}(\textbf{n})]}{n}= \lim_{n\rightarrow
\infty}\frac{\tilde{T}(\textbf{n})}{n},
\end{equation}

\vspace{0.1cm}

\noindent where $\textbf{n}=(n,0,\ldots,0)$. The existence of the
limit in (\ref{eq:mu}) and the fact that it is an almost sure
constant is proved in Deijfen (2003). The proof does not use the
assumption of bounded support for the radius distribution and
hence it applies also if this assumption is dropped. In fact, the
only part of the proof of the shape theorem in Deijfen (2003) that
uses the assumption of bounded support for $F$ is the one that
establishes that $\mu>0$, that is, that the infection does not
grow faster than linearly in time. Hence a weaker condition that
guarantees at most linear growth could replace the bounded support
assumption without weakening the conclusion of the theorem. We
will show that existence of the moment generating function of $F$
is sufficient for the growth to be at most linear. More precisely,
we will show:

\begin{prop} \label{prop:mugeq0}
If $\int_0^\infty e^{-\varphi r}dF(r)<\infty$ for some
$\varphi<0$, then $\mu>0$.
\end{prop}

\noindent In view of the above discussion, Theorem
\ref{th:aformgen} follows once we have proved Proposition
\ref{prop:mugeq0}. The main ingredient in the proof of Proposition
\ref{prop:mugeq0} is a ``larger'' growth process in which the
outbursts constitute a spatial branching process. This process
will be referred to as the Branching Random Walk growth process.
We will show that the BRW process grows at most linearly in time
and, since it can be shown that the original growth process is
stochastically dominated by the BRW process, the proposition
follows from this. The time constant $\mu$ is defined based on a
unit rate process and hence we consider only unit rate processes
for the remainder of this section.\medskip

\noindent The BRW growth process works in a similar way as the
original process, with outbursts that infect a randomly sized
shape around the outburst point. In the BRW process though, each
outburst point is assigned its own independent Poisson process to
generate new outbursts in the surrounding outburst shape.
Furthermore, for technical reasons we will take the outbursts in
the BRW process to be cubes rather than spheres.\medskip

\noindent To formally construct the BRW growth process, let
$\{N_n\}$ be a sequence of independent unit rate Poisson processes
on $\mathbb{R}^{d+1}$. The extra dimension represents time and the
points in the $n$th process are denoted $(X_k^n,T_k^n)$ where
$X_k^n\in\mathbb{R}^d$ and $T_k^n$ gives the location on the time
axis. Also, to each Poisson point, associate independently a
variable $R_k^n$ with distribution $F$. Finally, for $S\subset
{\mathbb{R}}^{d}$, let $N_n(S\times \mathbb{R})$ denote the
restriction of $N_n$ to $S\times {\mathbb{R}}$. The process now
evolves at time points $\{T_n\}$ by aid of cubic outbursts with
side length $\{2R_n\}$ centered at points $\{X_n\}$ obtained
inductively as follows:

\begin{itemize}
\item[1.] Define $X_0=0$, $T_0=0$ and $R_0=\gamma$ and let $C_n$ denote a cube in
${\mathbb{R}}^{d}$ with side length $2R_n$ centered at $X_n$.
\item[2.] Given $\{X_i; \hspace{0.1cm} i\leq n\}$,
$\{T_i; \hspace{0.1cm} i\leq n\}$ and $\{R_i; \hspace{0.1cm} i\leq
n\}$, for $i=0,\ldots ,n$, let
$$
\acute{T}_i^n=\inf_{k}\{T_k^i;\hspace{0.1cm} T_k^i>T_n \textrm{
and } (X_k^i,T_k^i)\in N_i(C_i\times {\mathbb{R}})\}
$$
and define $T_{n+1}=\min_i\{\acute{T}_i^n\}$. The point $X_{n+1}$
is the (a.s.\ unique) point in ${\mathbb{R}}^{d}$ such that
$(X_{n+1},T_{n+1})\in N_i$ for some $i$ and $R_{n+1}$ is the side
length variable associated with $(X_{n+1},T_{n+1})$.
\end{itemize}

\noindent The infected region after $n$ outbursts is obtained as
$\bar{S}_{(n)}=\cup_{i=0}^{n}{C_i}$ and for $t\in [T_n,T_{n+1})$
the infected region at time $t$ is given by
$\bar{S}_t=\bar{S}_{(n)}$.\medskip

\noindent \textbf{Remark 2.1} In the above construction the
initial set $\bar{S}_0$ is a cube with side length $2\gamma$
centered at the origin. The notation $\bar{S}_t$ is reserved for
the infected region at time $t$ starting from this particular
configuration. Furthermore, in the original one-type model, it
will often be convenient to take $S_0=B(0,\gamma)$ and the
notation $S_t$ is henceforth used to represent the infected region
at time $t$ starting from this particular choice of $S_0$.\medskip

\noindent The first result is a lemma stating that $\{S_t\}$ is
stochastically dominated by $\{\bar{S}_t\}$.

\begin{lemma}\label{lemma:brwdom}
The processes $S_t$ and $\bar{S}_t$ can be coupled in such a way
that $S_t\subset \bar{S}_t$ for all $t$.
\end{lemma}

\noindent \emph{Proof:} First note that, by definition, we have
$S_0=B(0,\gamma)$ and $\bar{S}_0=C(0,2\gamma)$, where
$C(0,2\gamma)$ is a cube with side length $2\gamma$ centered at
the origin. Hence $S_0\subset \bar{S}_0$. Let $N_0$ be the Poisson
process used to generate $\{S_t\}_{t>0}$ and let $\{N_k\}_{k\geq
1}$ be a sequence of independent Poisson processes on
$\mathbb{R}^{d+1}$ that are independent also of $N_0$. A process
with the same distribution as $\{\bar{S}_t\}$ is obtained by
starting from $C(0,2\gamma)$ and then using $\cup_{i=0}^{n-1}N_i$
to generate new outbursts in regions that previously have been
exposed to $n$ outbursts. That is, in intersections between $n$
outburst balls we scan within $n$ independent Poisson processes,
always including $N_0$, upwards along the time axis to find new
outburst points. Note that this is a different (but equivalent)
way of constructing the process compared to the above definition.
Some thought reveals that with this construction, given that
$S_{(n)}\subset\bar{S}_{(n)}$, it will also hold that
$S_{(n+1)}\subset\bar{S}_{(n+1)}$. It follows by induction over
$n$ that $S_{(n)}\subset\bar{S}_{(n)}$ for all $n$, and the lemma
is thereby proved. $\hfill\Box$\medskip

\noindent The outbursts in the BRW process satisfy the
independence structure usually assumed in branching processes, and
is in fact a branching process with no deaths of the
Crump-Mode-Jagers type; see e.g.\ Chapter 6 in Jagers (1975) for a
description of the general process. Note that the reproduction of
the ancestor at the origin is slightly different from the other
individuals reproduction in that $R_0\equiv \gamma$, that is, the
side length of the infection cube surrounding the ancestor is
deterministic.\medskip

\noindent Before proving Proposition \ref{prop:mugeq0} we state a
theorem by Biggins that will play a key role in the proof. To
formulate the result, consider a one-dimensional general spatial
branching process in which all individuals are equal, that is, all
individuals are of the same type and the distribution of its
progeny in space and time is the same. The reproduction of an
individual is described by a point process $Z$ on
$\mathbb{R}\times \mathbb{R}^+$ with each point corresponding to a
child. Let the intensity measure of $Z$ be denoted $\nu$ and let
$m(\varphi,\phi)$ be its Laplace transform, that is,
$$
m(\varphi,\phi)=\int e^{-\varphi x-\phi t}\nu(dx,dt).
$$
Define
$$
\alpha(\varphi)=\inf\{\phi;\, m(\varphi,\phi)\leq 1\}.
$$
\noindent Finally, write $H_t$ for the position of the rightmost
individual at time $t$.

\begin{theorem}[Biggins 1995] \label{th:Biggins} Assume that
$\alpha(\varphi)<\infty$ for some $\varphi<0$. Then there is a
constant $\zeta<\infty$ such that almost surely $H_t/t\rightarrow
\zeta$ as $t\rightarrow \infty$.
\end{theorem}

\noindent With this result at hand we are ready to prove
Proposition \ref{prop:mugeq0}.\medskip

\noindent \emph{Proof of Proposition \ref{prop:mugeq0}:} We want
to prove that the infected region in the original process grows at
most linearly in time. By Lemma \ref{lemma:brwdom} this follows if
we can show that the growth of the BRW process is at most linear.
To do this we will assume that the ancestor in the BRW growth
process has the same reproduction as the other individuals, that
is, we will assume that the process is started from a cube with
random side length distributed according to $F$. Linear growth for
such a process establishes linear growth also for a process with
an ancestor surrounded by a deterministic cube. This can be seen
as follows: Let $\bar{S}_t^\Gamma$ denote the infected region at
time $t$ in a BRW process started from an arbitrary initial set
$\Gamma$ and let $\bar{S}_t^{[\Gamma,s]}$ denote the region
infected at time $t\geq s$ in a BRW process started at time $s$
emanating from $\Gamma$. If $\tau$ denotes the time when the cube
$C(0,2\gamma)$ is infected in a BRW process started from a cube
with random side length $2R$, where $R\sim F$, then clearly
$$
\bar{S}^{[C(0,2\gamma),\tau]}_t\subset \bar{S}_t^{C(0,2R)}
$$
\noindent for $t\geq \tau$ and, since
$\bar{S}^{[C(0,2\gamma),\tau]}_t$ has the same distribution as
$\bar{S}_{t-\tau}$, linear growth for $\bar{S}_t^{C(0,2R)}$
guarantees linear growth also for $\bar{S}_t$.\medskip

\noindent To prove linear growth in a BRW growth process with
i.i.d.\ reproductions for all individuals, including the ancestor,
consider the projection of such a process on the first coordinate
axis. Some thought reveals that this projection is a
one-dimensional branching process in which each individual gives
birth to children according to a Poisson process in time with rate
$(2R)^d$, where $R$ is a random variable with distribution $F$.
The children are distributed uniformly in an interval of length
$2R$ centered at the parent. Given that $R=r$ in such a process,
we have $\nu(dx,dt)=(2r)^{d-1}dxdt$ and the Laplace transform
becomes

\begin{eqnarray*}
m_r(\varphi,\phi) & = & \int_0^\infty\int_0^{2r}
(2r)^{d-1}e^{-\varphi x-\phi t}dxdt\\[0.1cm]
& = & (\varphi\phi)^{-1}(2r)^{d-1}\left(1-e^{-2\varphi r}\right).
\end{eqnarray*}

\noindent Integrating over $r$ yields
$$
m(\varphi,\phi)=(\varphi\phi)^{-1}\int_0^\infty
(2r)^{d-1}\left(1-e^{-2\varphi r}\right)dF(r),
$$
which implies that
$$
\alpha(\varphi)=\varphi^{-1}\int_0^\infty
(2r)^{d-1}\left(1-e^{-2\varphi r}\right)dF(r).
$$
Hence $\alpha(\varphi)<\infty$ iff $\int_0^\infty e^{-\varphi
r}dF(r)<\infty$. Assume that this is the case for some $\varphi<0$
and let $H_t^1$ denote the position of the rightmost individual in
the projected process. Then, by Theorem \ref{th:Biggins}, there is
a constant $\zeta$ such that $H_t^1/t\rightarrow \zeta$ as
$t\rightarrow \infty$. Thus, if (\ref{eq:rvillk}) holds, the BRW
process grows at most linearly in time in the direction of the
first coordinate axis. But the same reasoning can be applied to
all coordinate axes: If $H_t^i$ denotes the position of the
rightmost individual in the projection of the BRW process on the
$i$:th coordinate axis ($i=1,\ldots ,d$), we have that
$H_t^i/t\rightarrow \zeta$ for each $i$ as $t\rightarrow \infty$.
Hence, for any $\varepsilon>0$, on the scale $1/t$ the set of
outbursts in the BRW growth process is contained in a cube with
side length $2\zeta+\varepsilon$ centered at the origin if $t$ is
sufficiently large. This means that the process grows at most
linearly in time and Lemma \ref{lemma:brwdom} completes the
proof.$\hfill\Box$

\section{Construction of the two-type model}

\noindent In this section the two-type model is built up more
formally by the construction of a Markov process whose state at
time $t$ is a subset of $\mathbb{R}^d$ and consists of two
disjoint sets $S_t^1$ and $S_t^2$. The process may for example be
thought of as describing the growth of two competing germ colonies
and the set $S_t^1$ ($S_t^2$) will be referred to as the type 1
(2) infected region.\medskip

\noindent To construct the model, let $N_1$ and $N_2$ be two
independent Poisson processes on $\mathbb{R}^{d+1}$ with
intensities $\lambda_1$ and $\lambda_2$ respectively. The extra
dimension represents time and the points in $N_i$ ($i=1,2$) are
denoted $(X_k^i,T_k^i)$, where $X_k^i\in \mathbb{R}^d$ and $T_k^i$
gives the location on the time axis. Furthermore, to each point in
the Poisson processes a random radius is associated. The radius
variables are assumed to be i.i.d.\ with distribution $F$ and
expected value $\gamma$. We will use the processes $N_1$ and $N_2$
together with the attached radius variables to construct three
sequences $\{T_n\}$, $\{X_n\}$ and $\{R_n\}$ indicating the time
points, locations and radii respectively of the outbursts, and two
sequences $\{S_{(n)}^i\}$ ($i=1,2$) specifying the type $i$
infected region after $n$ outbursts. The intuition is as follows:
At time zero a ball with radius $\gamma$ around the point
$-\textbf{2}\gamma= (-2\gamma,0,\ldots ,0)$ is infected with type
1 infection and a ball with radius $\gamma$ around the origin is
infected with type 2 infection so that
$S_{(0)}^1=B(-\textbf{2}\gamma,\gamma)$ and
$S_{(0)}^2=B(0,\gamma)$. The growth is then brought about by
scanning within the set $(S_{(0)}^1\cup S_{(0)}^2)\times
\mathbb{R}$ upwards along the time axis until either $S_{(0)}^1$
hits a point in $N_1$ or $S_{(0)}^2$ hits a point in $N_2$. An
outburst then takes place at this location, infecting all points
within some random distance from the outburst point. The type of
the infection is determined by the region in which the outburst
occurs: An outburst in the type $i$ infected region generates
outburst balls of type $i$. After the outburst the new infected
region is given by $S_{(1)}^1\cup S_{(1)}^2$, where the infected
region with the same infection type as the outburst might be
enlarged compared to before the outburst and the other region is
unchanged. Next we follow the set $(S_{(1)}^1\cup S_{(1)}^2)\times
\mathbb{R}$ further upwards along the time axis and eventually one
of the regions $S_{(1)}^1$ and $S_{(1)}^2$ hits a new point in
$N_1$ or $N_2$ respectively. This causes a new outburst and the
infected region is enlarged in the same way as described above.
And so on.\medskip

\noindent Formally the sequences $\{T_n\}$ (time points for the
outbursts), $\{X_n\}$ (locations of the outbursts), $\{R_n\}$
(radii of the outburst balls) and $\{S_{(n)}^i\}$ are constructed
inductively as follows:

\begin{itemize}
\item[1.] Define $T_0=0$, $S_{(0)}^1=B(-\textbf{2}\gamma,\gamma)$
and $S_{(0)}^2=B(0,\gamma)$. Also, for $S\subset \mathbb{R}^d$,
let $N_i(S\times \mathbb{R})$ denote the restriction of $N_i$ to
$S\times \mathbb{R}$.

\item[2.] Given $T_n$ and $S_{(n)}^i$ ($i=1,2$), define
$T_{n+1}=\min\{\acute{T}_{n+1}^1, \acute{T}_{n+1}^2\}$, where
$$
\acute{T}_{n+1}^i=\inf_{k}\{T_k^i; \hspace{0.1cm} T_k^i>T_n
\textrm{ and } (X_k^i,T_k^i)\in N_i(S_{(n)}^i\times \mathbb{R})\}.
$$
\noindent The point $X_{n+1}$ is the (a.s.\ unique) point in
$\mathbb{R}^d$ such that $(X_{n+1},T_{n+1})\in N_i$ for some $i$
and $R_{n+1}$ is the radius variable associated with the point
$(X_{n+1},T_{n+1})$.

\item[3.] Once the points $T_{n+1}$, $X_{n+1}$ and $R_{n+1}$ are
specified, the infected regions $S_{(n)}^1$ and $S_{(n)}^2$ are
updated as follows: If $(X_{n+1},T_{n+1})\in N_1$, that is, if the
outburst is of type 1, then
$$
\left\{
\begin{array}{l}
S_{(n+1)}^1=S_{(n)}^1\cup[B(X_{n+1},R_{n+1})\cap (S_{(n)}^1\cup
S_{(n)}^2)^c]\\[0.1cm]
S_{(n+1)}^2=S_{(n)}^2
\end{array}
\right.
$$

\noindent so that, in words, the type 1 infected region is
enlarged by the previously uninfected parts of the outburst ball
$B(X_{n+1},R_{n+1})$ and the type 2 infected region remains
unchanged. If the outburst is of type 2, that is, if
$(X_{n+1},T_{n+1})\in N_2$, then the type 2 infected region is
updated analogously while the type 1 infected region is left
unchanged.
\end{itemize}

\noindent For $t\in [T_n,T_{n+1})$, the type $i$ infected region
at time $t$ is given by $S_t^i=S_{(n)}^i$ and the total infected
region at time $t$ is $S_t^1\cup S_t^2$. \medskip

\noindent Write $\Delta_{n+1}^i$ for the time counting from $T_n$
until an outburst of type $i$ occurs. By standard properties of
the Poisson process, we have
$$
\Delta_{n+1}^i|{\mathcal{F}}_{n}^i \sim
\textrm{Exp}\{\lambda_i|S_{(n)}^i|\},
$$
\noindent where ${\mathcal{F}}_n^i=\sigma (S_{(0)}^i,\ldots
,S_{(n)}^i)$. Furthermore, given $S_t^i$ the memoryless property
of the exponential distribution implies that the time until an
outburst occurs somewhere in $S_t^i$ is exponentially distributed
with parameter $\lambda_i|S_t^i|$ and, also, the location of the
outburst is uniformly distributed over $S_t^i$.\medskip

\noindent To make sure that the model does not explode by
generating infinitely many outbursts in finite time, we have to
show that the sequence $\{T_n\}$ does not have a finite limit
point. This is done in the following proposition, which is the
analogue of Proposition 2.1 in Deijfen (2003).

\begin{prop}\label{prop:ejexpl}
If the radius distribution has finite moment of order $d$, then
almost surely $T_n\rightarrow \infty$ as $n\rightarrow \infty$.
\end{prop}

\noindent \emph{Proof:} Let $\{\Delta_n\}$ denote the increments
of the process $\{T_n\}$, that is, $\Delta_n:=T_{n}-T_{n-1}$ is
the time between two successive outbursts regardless of type.
Since $T_n=\sum_{k=1}^n\Delta_k$, the proposition follows if we
can show that $\sum_{k=1}^\infty \Delta_k=\infty$ almost surely.
To do this, note that
$$
\Delta_k=\min\{\Delta_k^1,\Delta_k^2\}.
$$
\noindent Since $\Delta_{k}^i|{\mathcal{F}}_{k-1}^i \sim
\textrm{Exp}\{\lambda_i|S_{(k-1)}^i|\}$, it follows that
$$
\Delta_k|{\mathcal{F}}_{k-1}\sim
\textrm{Exp}\{\lambda_1|S_{(k-1)}^1|+\lambda_2|S_{(k-1)}^2|\},
$$
\noindent where ${\mathcal{F}}_{k}=\sigma(S_{(0)}^1,\ldots
S_{(k)}^1,S_{(0)}^2,\ldots S_{(k)}^2)$. Furthermore, by properties
of the Poisson process, given ${\mathcal{F}}_{k-1}$ we can write
$$
\Delta_k=\frac{k}{\lambda_1|S_{(k-1)}^1|+\lambda_2|S_{(k-1)}^2|}\cdot
E_k
$$
where $\{E_k\}$ are independent, $E_k\sim\textrm{Exp}(k)$. A
trivial upper bound for $|S_{(k-1)}^i|$ ($i=1,2$) is given by
$$
|S_{(k-1)}^i|\leq v_0+\sum_{n=1}^{k-1}V_n,
$$
\noindent where $v_0$ is the volume of the initial type $i$
$\gamma$-ball in $\mathbb{R}^d$ and $V_n$ denotes the volume of a
$d$-dimensional ball with radius $R_n$. Let $v=\textrm{E}[V_n]$,
which is finite by the assumption of the proposition. By the
strong law of large numbers,
$$
\frac{1}{k}\sum_{n=1}^{k-1}V_n\rightarrow v \quad \textrm{as } k
\rightarrow \infty
$$
\noindent and hence, if $k$ is large,
$$
\frac{1}{k}|S_{(k-1)}^i|\leq 2v.
$$
\noindent Thus, for large $k$,
$$
\Delta_k\geq \frac{1}{(\lambda_1+\lambda_2)2v}\cdot E_k
$$
\noindent and we are done if we can show that $\sum_{k=1}^\infty
E_k=\infty$. But this is an easy consequence of Kolmogorov's three
series theorem: Let $\tilde{E}_k=E_k-\textrm{E}[E_k]=E_k-1/k$, so
that $E_k=\tilde{E}_k+1/k$. Since $\sum_{k=1}^{\infty}1/k=\infty$
it suffices to show that $\sum_{k=1}^{\infty} \tilde{E}_k$
converges almost surely. Using the fact that
$\sum_{k=1}^{\infty}\textrm{E}[{\tilde{E}_k}^2]=
\sum_{k=1}^{\infty}1/k^2<\infty$ this follows from the three
series theorem. $\hfill \Box$

\section{Preliminaries}

\noindent In this section we prove a number of auxiliary results
needed in the later sections. The first two lemmas concern the
relation between the one-type model and the two-type model.

\begin{lemma}\label{lemma:svagdom} Let $S_t^2$ denote the type 2
infected region at time $t$ in a two-type process with
distribution $P^{1,\lambda}$ and let $S_t$ denote the infected
region at time $t$ in a one-type process with rate $\lambda$. The
two-type process and the one-type process can be coupled in such a
way that $S_t^2\subset S_t$ for all $t$.\end{lemma}

\noindent \emph{Proof:} Couple the two processes by letting the
one-type process and the type 2 outbursts in the two-type process
be generated by the same rate $\lambda$ Poisson process with the
same radius variables attached. If $S_{(n)}^2\subset S_{(n)}$, it
then also holds that $S_{(n+1)}^2\subset S_{(n+1)}$. Since
$S_0^2=S_0$, the lemma follows by induction over $n$.$\hfill
\Box$\medskip

\begin{lemma}\label{lemma:ttpetp}
Consider the one-type process $\{S_t\}_{t\geq 0}$ with rate
$\lambda\leq 1$ and the two-type process $\{S_t^1\cup
S_t^2\}_{t\geq 0}$ with distribution $P^{1,\lambda}$. These can be
coupled in such a way that

\begin{equation}\label{eq:l1stokdom}
S_t\subset S_t^1\cup S_t^2
\end{equation}

\noindent for all $t$.
\end{lemma}

\noindent \emph{Proof:} For $t=0$, (\ref{eq:l1stokdom}) is
trivial. To couple the one-type and the two-type processes, let
$N_1$ and $N_2$ be two independent Poisson processes with rate
$1-\lambda$ and $\lambda$ respectively. Use $N_1\cup N_2$ to
generate the type 1 outbursts in the two-type process and use
$N_2$ to generate all outbursts in the one-type process and the
type 2 outbursts in the two-type process. Note that the two-type
process is obtained in a different way here compared to in Section
3, but it is easy to see that it gets the correct distribution.
Also, it is easy to see that (\ref{eq:l1stokdom}) is preserved for
$t>0$. $\hfill\Box$\medskip

\noindent The next two lemmas are needed in the proof of
Proposition \ref{prop:key}. To formulate them, introduce a new,
hampered version of the one-type process by placing ``ceilings''
and ``floors'' in $\mathbb{R}^d$ restricting the growth in all
directions but one: Write $(x_1,\ldots ,x_d)$ for the coordinates
of a point $x\in\mathbb{R}^d$ and let $S_t^b$ denote the infected
region at time $t$ in a one-type process where all infection
outside the stripe $\Omega_b=\{x\in\mathbb{R}^d;\hspace{0.1cm}
|x_i|\leq b \textrm{ for all } i\geq 2\}$ is ignored. The process
$\{S_t^b\}$ hence works exactly like the original one-type process
except that points with $|x_i|\geq b$ for some $i\geq 2$ are
immune to the infection. The following lemma says that $\Omega_b$
is filled with infection linearly in time.

\begin{lemma}\label{lemma:formbeg} Consider a hampered one-type process
with unit rate. Assume that (\ref{eq:rvillk}) holds and let
$S_0^b\subset \Omega_b$ be bounded with strictly positive Lebesgue
measure. Then, for any dimension $d$, there is a real number
$\mu_b>0$ such that, for any $\varepsilon\in(0,\mu_b^{-1})$,
almost surely
$$
(1-\varepsilon)\left\{x\in\Omega_b;\,|x_1|\leq t
\mu_b^{-1}\right\}\subset S_t^b\subset
(1+\varepsilon)\left\{x\in\Omega_b;\,|x_1|\leq t
\mu_b^{-1}\right\}
$$
\noindent for all sufficiently large $t$.
\end{lemma}

\noindent The proof of the lemma for the case of bounded outburst
radii is a straightforward but tedious adaptation of the proof of
the shape theorem in Deijfen (2003), and the general case follows
as in Section 2. We therefore omit the proof.\medskip

\noindent Let $\tilde{T}^b(x)$ be the analogue of $\tilde{T}(x)$
in the process $S_t^b$, that is, $\tilde{T}^b(x)$ is the time when
the $\gamma$-ball around the point $x$ is infected in a unit rate
hampered process started from $S_0=B(0,\gamma)$. The time-constant
$\mu_b$ is defined analogously to the time-constant for the
unhampered process, that is,
$$
\mu_b:=\lim_{n\rightarrow \infty}
\frac{\textrm{E}[\tilde{T}^b(\textbf{n})]}{n}=\lim_{n\rightarrow
\infty}\frac{\tilde{T}^b(\textbf{n})}{n},
$$
\noindent where $\textbf{n}=(n,0,\ldots,0)$. The following lemma
states that, as $b$ becomes large, the speed of the growth in the
hampered process approaches the speed in the unhampered process.

\begin{lemma}\label{lemma:mubtmu}
As $b\rightarrow\infty$ we have $\mu_b\rightarrow\mu$.
\end{lemma}

\noindent \emph{Proof:} Trivially $\mu_b\geq \mu$ so it suffices
to show that $\lim_{b\rightarrow \infty}\mu_b\leq \mu$. To this
end, consider a one-type process with unit rate and pick
$\delta>0$ and $p\in(0,1)$. We will show that

\begin{equation}\label{eq:tubnk}
P\left(\tilde{T}^b(k\textbf{n})>(1+\delta)\mu nk\right)\leq p
\end{equation}

\noindent if $n$, $k$ and $b$ are large. Since $p>0$ was arbitrary
this implies that almost surely
$$
\lim_{n\rightarrow\infty} \frac{\tilde{T}^b(\textbf{n})}{n}\leq
(1+\delta)\mu
$$
for large $b$ and, since also $\delta>0$ was arbitrary, the
proposition follows. To prove (\ref{eq:tubnk}), first note that by
Theorem \ref{th:aformgen} and (\ref{eq:mu}) we have

\begin{equation}\label{eq:ttbegr}
\textrm{E}[\tilde{T}(\textbf{n})]\leq (1+\delta/3)\mu n,
\end{equation}

\noindent if $n$ is large. Define $D_n^b=\tilde{T}^b(\textbf{n})-
\tilde{T}(\textbf{n})$ and let $F^b_n$ be the event that the
hampered process $S_t^b$ reaches $\partial\Omega_b$ before time
$\tilde{T}^b(\textbf{n})$. We will show that

\begin{itemize}
\item[(i)] $P(F^b_n)\rightarrow 0$ as $b\rightarrow \infty$;
\item[(ii)] $\textrm{E}[D_n^b|F^b_n]\leq cn$ for some constant $c\in\mathbb{R}$.
\end{itemize}

\noindent The claim (i) follows easily by noting that
$P(F^b_n)\leq P(\|S_{\tilde{T}(\textbf{n})}\|>b)$. Since almost
surely $\tilde{T}(\textbf{n})<\infty$, Proposition
\ref{prop:ejexpl} gives that
$P(\|S_{\tilde{T}(\textbf{n})}\|<\infty)=1$ and hence
$P(\|S_{\tilde{T}(\textbf{n})}\|>b)\rightarrow 0$ as $b\rightarrow
\infty$.\medskip

\noindent To establish (ii), write $\tau_b$ for the time when the
infection reaches $\partial\Omega_b$ and note that
$$
\textrm{E}[D_n^b|F^b_n]\leq
\textrm{E}\left[\tilde{T}^b(\textbf{n})-\tau_b|F^b_n\right].
$$
\noindent Now imagine that at time $\tau_b$ a new process is
started from the origin using only outbursts that touch the
$x$-axis, that is, at time $\tau_b$ all infection except a ball
with radius $\gamma$ around the origin is erased and the infection
then evolves in time along the $x$-axis using the same
$d+1$-dimensional Poisson process as the original process. Let
$\tau_{\textbf{n}}$ denote the time, counting from $\tau_b$, when
the $\gamma$-ball around the point $\textbf{n}$ is infected in
this new process. Since $\tilde{T}^b(\textbf{n})\leq
\tau_b+\tau_{\textbf{n}}$ we have

\begin{equation}\label{eq:forsbegr}
\textrm{E}[D_n^b|F^b_n]\leq \textrm{E}[\tau_{\textbf{n}}].
\end{equation}

\noindent Using the same technique as in the proof of Lemma 3.1 in
Deijfen (2003), it follows that the time until the point
$\textbf{n}$ is infected in the $x$-axis process can be bounded by
a sum of $n\lceil 2\gamma^{-1}\rceil$ independent exponential
variables with mean $\eta=\eta(d)$. Furthermore, it is not hard to
see that the time from when the point $\textbf{n}$ is infected
until the entire $\gamma$-ball around $\textbf{n}$ is infected can
be bounded by a sum of $m=m(d)$ exponential variables, which may
be defined so that their mean equals $\eta$. Hence
$\textrm{E}[\tau_{\textbf{n}}]\leq cn$, where $c$ can be taken as
$\eta\lceil m+2\gamma^{-1}\rceil$. The statement (ii) now follows
from (\ref{eq:forsbegr}).\medskip

\noindent By (i) we can pick $b$ large so that $P(F^b_n)\leq
\mu\delta/3c$. Using (ii) and the fact that
$$
\textrm{E}[D_n^b]= P(F^b_n)\textrm{E}[D_n^b|F^b_n],
$$
\noindent it follows that, for such $b$, we have

\begin{equation}\label{eq:forvfors}
\textrm{E}[D_n^b]\leq \delta\mu n/3.
\end{equation}

\noindent Now, if $n$ is chosen large enough to ensure
(\ref{eq:ttbegr}) and $b$ large enough to ensure
(\ref{eq:forvfors}), then

\begin{eqnarray}
\textrm{E}[\tilde{T}^b(\textbf{n})] & = &
\textrm{E}[\tilde{T}(\textbf{n})+D_n^b]\nonumber \\
 & \leq & (1+\delta/3)\mu n+\delta\mu n/3\nonumber \\
 & = & (1+2\delta/3)\mu n.\label{eq:ettbn}
\end{eqnarray}

\noindent It remains to show that this implies (\ref{eq:tubnk}).
To this end, let $\tilde{T}^b((j-1)\textbf{n}, j\textbf{n})$
denote the time it takes for the infection to invade the entire
$\gamma$-ball around the point $j\textbf{n}$ in a process started
at time $\tilde{T}^b((j-1)\textbf{n})$ emanating from the point
$(j-1)\textbf{n}$. The variables
$\{\tilde{T}^b((j-1)\textbf{n},j\textbf{n});\hspace{0.1cm}j\geq
1\}$ are i.i.d.\ with expected value
$\textrm{E}[\tilde{T}^b(\textbf{n})]$ and hence, by the strong law
of large numbers, almost surely
$$
\frac{1}{k}\sum_{j=1}^{k}\tilde{T}^b((j-1)\textbf{n},j\textbf{n})\rightarrow
\textrm{E}[\tilde{T}^b(\textbf{n})]\quad \textrm{as } k\rightarrow
\infty.
$$
\noindent Using (\ref{eq:ettbn}) this implies that
$$
P\left(\sum_{j=1}^{k}\tilde{T}^b((j-1)\textbf{n},j\textbf{n})>
(1+\delta)\mu nk \right)\leq p
$$
\noindent if $k$ is large. Since
$$
\tilde{T}^b(k\textbf{n}) \leq
\sum_{j=1}^{k}\tilde{T}^b((j-1)\textbf{n},j\textbf{n}),
$$
this proves (\ref{eq:tubnk}). $\hfill \Box$\medskip

\noindent The next lemma is needed to prove Proposition
\ref{prop:startomr}. It involves the concept of \emph{effective}
outbursts: An outburst is said to be effective if it causes
previously uninfected regions to be infected, that is, if it
reaches outside the boundary of the infected region.

\begin{lemma}\label{lemma:andlutbr} Assume that $F$
satisfies (\ref{eq:rvillk}) and let $\Lambda$ be a bounded subset
of $\mathbb{R}^d$.

\begin{itemize}
\item[\rm{(a)}] The number of effective outbursts that occur in $\Lambda$
during the progress of the growth in a two-type process is almost
surely finite.
\item[\rm{(b)}] If in addition $F$ has unbounded support and if
$\Lambda^c\cap [S_0^1\cup S_0^2]$ has positive Lebesgue measure,
there is a positive probability that no effective outbursts ever
occur in $\Lambda$.
\end{itemize}
\end{lemma}

\noindent \textbf{Remark 4.1} In analogy with Remark 1.1, Lemma
\ref{lemma:andlutbr}(b) extends to the case with bounded radii
provided that $\Lambda^c \cap [S^1_0 \cup S^2_0]$ is not
``strangled'' by $\Lambda \cap [S^1_0 \cup S^2_0]$ in the sense of
Remark 1.1.\medskip

\noindent \emph{Proof of Lemma \ref{lemma:andlutbr}:} By
time-scaling and symmetry it is enough to consider a process with
infection rates $(1,\lambda)$, where $\lambda\leq 1$. Furthermore,
the choice of initial sets does not affect the arguments in the
proof. Hence we may restrict our attention to a process with
distribution $P^{1,\lambda}$.\medskip

\noindent (a) Write $N_\Lambda$ for the number of effective
outbursts in $\Lambda$. Lemma \ref{lemma:ttpetp} and Theorem
\ref{th:aformgen} implies that almost surely
$$
S_t^1\cup S_t^2\supset B\left(0,\frac{\lambda\mu^{-1}t}{2}\right)
$$
for all sufficiently large $t$. Hence the minimal distance between
points in $\Lambda$ and points in $(S_t^1\cup S_t^2)^c$ is at
least $\lambda\mu^{-1}t/4$ for all sufficiently large $t$. It
follows that

\begin{eqnarray*}
\textrm{E}[N_\Lambda] & \leq & |\Lambda|\int_0^\infty\left(
\int_{\lambda\mu^{-1}t/4}^\infty dF(r)\right)dt\\
& = & |\Lambda|\int_0^\infty\left(\int_0^{4r/\lambda\mu^{-1}}dt\right)dF(r)\\
& = & \frac{4|\Lambda|}{\lambda\mu^{-1}}\int_0^\infty rdF(r).
\end{eqnarray*}

\noindent The last integral is finite by the assumption on $F$ and
hence $N_\Lambda$ is finite almost surely.\medskip

\noindent (b) The calculation in (a) shows that there is an $r$
(depending on $\Lambda$) such that, if $S_t^1\cup S_t^2$ contains
the ball $B(0,r)$, then the conditional expectation of the number
of effective outbursts in $\Lambda$ after time $t$ is at most
$1/2$. Let $A$ denote the event that the ball $B(0,r)$ is fully
infected before the first outburst in $\Lambda$. Then

\begin{eqnarray*}
P(N_\Lambda=0) & \geq & P(N_\Lambda=0|A)P(A)\\
& \geq & \frac{1}{2}P(A) ,
\end{eqnarray*}

\noindent which is clearly positive. \hfill $\Box$\medskip

\noindent We will later on need the following refinement of Lemma
\ref{lemma:andlutbr}.

\begin{lemma}\label{lemma:andlutbrqu} Consider a two-type process with
$F$ satisfying (\ref{eq:rvillk}). For any $\delta,\xi>0$, there
exists an $r^*_0<\infty$ such that for all $r^*\geq r^*_0$ we get
that, if the process is started with $(S_0^1,S_0^2)$ satisfying
$$
S_0^2\subset B(0,r^*)\quad \textrm{and}\quad S_0^1\cup
S_0^2\supset B(0,r^*(1+\xi)),
$$
then
$$
P(\textrm{infection 2 never makes an effective
outburst})>1-\delta.
$$
\end{lemma}

\noindent To prove Lemma \ref{lemma:andlutbrqu} we need the
following auxiliary result, which asserts that, if the initial set
in a two-type process is large, then the infection will continue
to grow at least with the speed stipulated by the shape theorem
for the weaker infection type.

\begin{lemma}\label{lemma:andlutbrquhjalp} For $\lambda\leq 1$ and
any $\delta,\varepsilon\in (0,1)$, there is an $s<\infty$ such
that, if a two-type process with infection rates $(1,\lambda)$ is
started in such a way that $S_0^1\cup S_0^2\supset
B(0,s\lambda\mu^{-1})$, then
$$
P\left(S_t^1\cup S_t^2\supset
B\left(0,(1-\varepsilon)(s+t)\lambda\mu^{-1}\right)\hspace{0,2cm}
\forall t\geq 0\right)>1-\delta.
$$
\end{lemma}

\noindent \emph{Proof of Lemma \ref{lemma:andlutbrquhjalp}:} In
view of Lemma \ref{lemma:ttpetp} it is enough to prove the
corresponding  statement for the one-type process, that is, it is
enough to prove that there is an $s<\infty$ such that, if a
one-type process with parameter $\lambda$ is started with
$S_0\supset B(0,s\lambda\mu^{-1})$, then

\begin{equation}\label{eq:hl1}
P\left(S_t\supset
B\left(0,(1-\varepsilon)(s+t)\lambda\mu^{-1}\right)\hspace{0,2cm}\forall
t\geq 0\right)>1-\delta.
\end{equation}

\noindent To do this, consider a one-type process
$\{S_t^*\}_{t\geq 0}$ with parameter $\lambda$ and initial
condition, say, $S_0^*=B(0,\gamma)$. Define the event
$$
A_t=\left\{
B\Big(0,\Big(1-\frac{\varepsilon}{2}\Big)t\lambda\mu^{-1}\Big)
\subset S_t^* \subset
B\Big(0,\Big(1+\frac{\varepsilon}{2}\Big)t\lambda\mu^{-1}\Big)
\right\},
$$
and note that, by Theorem \ref{th:aformgen}, there exists an
$s<\infty$ such that

\begin{equation}\label{eq:hl2}
P\bigg(A_t\textrm{ holds for all }t\geq
s\Big(1-\frac{\varepsilon}{2}\Big)\bigg)>1-\delta.
\end{equation}

\noindent Now couple the processes $\{S_t\}_{t\geq 0}$ and
$\{S_t^*\}_{t \geq 0}$ in such a way that the latter is generated
by the former's underlying Poisson process delayed by time
$s(1-\varepsilon/2)$. If $S_0\supset B(0,s\lambda\mu^{-1})$, on
the event in (\ref{eq:hl2}) we get that $S_0\supset
S^*_{s(1-\varepsilon/2)}$ and, by the choice of the coupling,
$S_t\supset S^*_{s(1-\varepsilon/2)+t}$ for all $t\geq 0$. This
implies (\ref{eq:hl1}).\hfill $\Box$\medskip

\noindent \emph{Proof of Lemma \ref{lemma:andlutbrqu}:} By
time-scaling and symmetry it suffices to consider a process with
infection rates $(1,\lambda)$, where $\lambda\leq 1$. If such a
process is started with $S_0^1\cup S_0^2\supset B(0,r^*(1+\xi))$,
it follows from Lemma \ref{lemma:andlutbrquhjalp} (with
$\varepsilon=\min\{\frac{\xi}{4},\frac{1}{4}\}$) that, if $r^*$ is
taken to be sufficiently large, then

\begin{equation}\label{eq:q1}
P\left(S_t^1\cup S_t^2\supset
B\Big(0,r^*\Big(1+\frac{\xi}{2}\Big)+(1-\varepsilon)t\lambda\mu^{-1}\Big)
\hspace{0,15cm} \forall t\geq 0\right)>1-\frac{\delta}{2}.
\end{equation}

\noindent Let $\pi_d$ be such that $\pi_d(r^*)^d$ is the volume of
a $d$-dimensional ball with radius $r^*$ and write $N_{r^*}$ for
the number of effective outbursts in $B(0,r^*)$. On the event in
(\ref{eq:q1}), we have that an effective outburst inside
$B(0,r^*)$ at time $t$ has to have radius at least
$r^*\xi/2+(1-\varepsilon)t\lambda\mu^{-1}$. Hence, on the event in
(\ref{eq:q1}) we have

\begin{eqnarray*}
\textrm{E}[N_{r^*}] & \leq & \pi_d(r^*)^d\int_0^\infty\left(
\int_{r^*\xi/2+(1-\varepsilon)t\lambda\mu^{-1}}^\infty dF(r)\right)dt\\
& \leq &
\pi_d(r^*)^d\int_{r^*\xi/2}^\infty\left(\int_0^{r/(1-\varepsilon)\lambda
\mu^{-1}}dt\right)dF(r)\\
& = & \frac{\pi_d(r^*)^d}{(1-\varepsilon)\lambda\mu^{-1}}
\int_{r^*\xi/2}^\infty rdF(r).
\end{eqnarray*}

\noindent When $r\geq r^*\xi/2$, we have $(r^*)^d\leq
2^dr^d/\xi^d$. Thus
$$
\textrm{E}[N_{r^*}]\leq
\frac{\pi_d2^d}{(1-\varepsilon)\lambda\mu^{-1}\xi^d}
\int_{r^*\xi/2}^\infty r^{d+1}dF(r).
$$
The last integral is finite by the assumption on $F$ and hence
E$[N_{r^*}]$ can be made arbitrarily small by taking $r^*$ large.
Take $r^*$ large enough so that E$[N_{r^*}]$ is at most $\delta/2$
and such that (\ref{eq:q1}) holds. Then, since $S_0^2\subset
B(0,r^*)$, with probability at least $1-\delta$ the type 2
infection never makes an effective outburst, as desired.\hfill
$\Box$\medskip

\section{A key proposition}

\noindent In this section we formulate and prove an analogue of
Proposition 2.2 in H\"aggstr\"om and Pemantle (2000). The
proposition will play a key role in the proof of Theorem
\ref{th:huvudres} and as in H\"aggstr\"om and Pemantle (2000), the
proof is rather lengthy and technical; this appears to be
unavoidable. We will opt for a geometrical argument that is a bit
different from the one of H\"aggstr\"om and Pemantle. The
unboundedness of the outbursts radii causes some extra
complications in our case, but on the other hand the fact that the
asymptotic shape is a sphere makes the geometric intuition a bit
more accessible.\medskip

\noindent We begin by observing that the events $G_1$ and $G_2$
have positive probability.

\begin{prop}\label{prop:pgisn} If $F$ has unbounded
support and satisfies (\ref{eq:rvillk}), then, for all infection
rates $(\lambda_1,\lambda_2)$ and all choices of initial sets
$(\Gamma_1,\Gamma_2)$ which are bounded and have positive Lebesgue
measure, we have
$P_{\Gamma_1,\Gamma_2}^{\lambda_1,\lambda_2}(G_i)>0$ for $i=1,2$.
\end{prop}

\noindent \textbf{Remark 5.1} Proposition \ref{prop:pgisn} extends
to the case of bounded radii as well, provided that neither
$\Gamma_1$ nor $\Gamma_2$ surrounds the other in the sense of
Remark 1.1.\medskip

\noindent \emph{Proof of Proposition \ref{prop:pgisn}:} Since the
total infected region increases to cover all of $\mathbb{R}^d$, we
have $P_{\Gamma_1,\Gamma_2}^{\lambda_1,\lambda_2}(G_1\cup G_2)=1$.
Furthermore, by Lemma \ref{lemma:andlutbr}(b), there is a positive
probability that no effective outbursts ever occur in $S_0^1$.
Hence $P_{\Gamma_1,\Gamma_2}^{\lambda_1,\lambda_2}(G_2)>0$.
Similarly it can be seen that
$P_{\Gamma_1,\Gamma_2}^{\lambda_1,\lambda_2}(G_1)>0$. \hfill
$\Box$\medskip

\noindent The proof of Theorem \ref{th:huvudres} is based on the
fact that, if both infection types will survive in the long run,
they have to grow equally fast, that is, the asymptotic speed of
the growth for the type 1 and the type 2 infection have to be the
same. This is formulated in Lemma \ref{lemma:keylemma}, which says
that on the event of indefinite growth for the weaker infection
type, the size of the asymptotic shape of the total infected
region is determined by the weaker infection type. The key step in
proving Lemma \ref{lemma:keylemma} is to show that, if the
stronger infection type gets a large enough lead over the weaker
infection type infinitely often -- which indeed will be the case
if the asymptotic shape is larger than the capacity of the weaker
infection type allows for -- then almost surely the stronger
infection type will eventually eradicate the weaker one. To
formulate this key result, assume that $\lambda_1=1$ (note that by
time-scaling this is no restriction) and write
$\lambda_2=\lambda$. Furthermore, for an arbitrary set
$\Gamma\subset \mathbb{R}^d$, let
$$
\|\Gamma\|=\sup\{|x|;\hspace{0.1cm} x\in \Gamma\}.
$$
The proposition now runs as follows.

\begin{prop}\label{prop:key} Assume that $F$ satisfies (\ref{eq:rvillk}).
Also, fix $\lambda<1$ and $\varepsilon >0$ and let
$$
D=\left\{\|S_t^1\|\geq (1+3\varepsilon)t\lambda\mu^{-1}\textrm{
for arbitrarily large }t\right\}.
$$
Then $P^{1,\lambda}(G_2|D)=0$.
\end{prop}

\noindent \emph{Proof:} Write $\tilde{S}_t^i$ for the set of
points whose entire $\gamma$-ball is contained in the type $i$
infected area at time $t$. Points in $\tilde{S}_t^i$ will be
referred to as \emph{strongly} type $i$ infected at time $t$. Let
$$
Q=(1+3\varepsilon)B(0,\lambda\mu^{-1})\backslash
(1+2\varepsilon)B(0,\lambda\mu^{-1})
$$
\noindent and introduce the event
$$
E_t=\{\tilde{S}_t^1\cap tQ\neq\emptyset\}.
$$
Note that almost surely

\begin{equation}\label{eq:dimplstark}
D\Rightarrow \{E_t\textrm{ occurs for arbitrarily large } t\}.
\end{equation}

\noindent To see this, write $\Pi$ for the set of type 1 outbursts
that occur in $tQ$ for some $t$ during the progress of the growth
and let $\Pi_\gamma$ be those outbursts in $\Pi$ whose radius is
at least $\gamma$. If $D$ occurs, then $|\Pi|=\infty$ and, since
each outburst in $\Pi$ has radius greater than $\gamma$ with some
probability $p>0$, it follows from Levy's version of the
Borel-Cantelli lemma (see Williams (1991), section 12.15) that
$|\Pi_\gamma|=\infty$ as well. But if $|\Pi_\gamma|=\infty$, the
region $tQ$ must contain strongly type 1 infected points
infinitely often and (\ref{eq:dimplstark}) is verified.\medskip

\noindent Now fix $\varepsilon >0$ and $\lambda< 1$. We want to
pick $\delta>0$ and $\alpha\in(0,\varepsilon]$ so that
$(1+\delta)^{-1}\mu^{-1}>(1+\alpha)\lambda \mu^{-1}$ (here the
left-hand side should be thought of as a lower bound for the speed
of a hampered unit rate process and the right-hand side as an
upper bound for the speed of an unhampered process with rate
$\lambda$). Hence we define
$$
\delta=\frac{1-\lambda}{2\lambda}
$$
and
$$
\alpha=\min\left\{\frac{1-\lambda}{2(1+\lambda)},\varepsilon\right\}.
$$
Let
$$
F_t=\{S_s^2\subset (1+\alpha)sB(0,\lambda\mu^{-1})\textrm{ for all
} s\geq t \},
$$
and write $\mathcal{F}_t=\sigma(S_s^1\cup
S_s^2;\hspace{0.1cm}s\leq t)$. We will show that, for some fixed
$c>0$, we have almost surely on the event $E_t$ that

\begin{equation}\label{eq:kritbeg}
P^{1,\lambda}(G_2^c|\mathcal{F}_t,F_t)>c
\end{equation}

\noindent if $t$ is large. Using Levy's 0-1 law the proposition
follows from this: By Theorem \ref{th:aformgen} and Lemma
\ref{lemma:svagdom}, $P^{1,\lambda}(F_t)\rightarrow 1$ as
$t\rightarrow \infty$. Together with (\ref{eq:dimplstark}) and
(\ref{eq:kritbeg}) this implies that, almost surely on the event
$D$,

\begin{equation}\label{eq:mennuda}
P^{1,\lambda}(G_2^c|\mathcal{F}_t)>c/2\quad \textrm{infinitely
often}.
\end{equation}

\noindent Levy's 0-1 law tells us that almost surely
$P^{1,\lambda}(G_2^c|\mathcal{F}_t)$ tends to the indicator
function of $G_2^c$ and (\ref{eq:mennuda}) prevents
$P^{1,\lambda}(G_2^c|\mathcal{F}_t)$ from converging to 0 on $D$.
Hence $P^{1,\lambda}(G_2^c|\mathcal{F}_t)\rightarrow 1$ on $D$,
implying that $P^{1,\lambda}(G_2|D)=0$, as desired.\medskip

\noindent It remains to prove (\ref{eq:kritbeg}). To understand
the idea of the proof, note that on $E_tF_t$ we have
$\|S_t^2\|\leq (1+\varepsilon)t\lambda\mu^{-1}$ and
$\|\tilde{S}_t^1\|\geq (1+2\varepsilon)t\lambda\mu^{-1}$, that is,
the strongly type 1 infected region at time $t$ has a lead of at
least $\varepsilon t\lambda \mu^{-1}$ units of length as compared
to the type 2 infected region. We will show that if $t$ is large,
then with large probability this lead gives the type 1 infection
time to create a layer that completely surrounds the type 2
infection. Moreover, if this layer is sufficiently thick -- which
it will indeed be if $t$ is large -- then Lemma
\ref{lemma:andlutbrqu} gives a lower bound for the probability
that no type 2 outbursts that reach outside the layer ever occur.
The proof is to a large extent based on a geometrical
construction, which is easiest to picture in two dimensions. Hence
we give the details for $d=2$ and indicate at the end of the proof
how the geometrical arguments can be generalized to higher
dimensions.\medskip

\begin{figure}
\centering \mbox{\epsfig{file=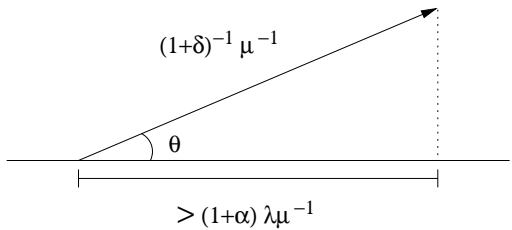,height=3cm}}
\caption{Choice of $\theta$.}
\end{figure}

\noindent To describe the geometrical construction, first define
an angle $\theta\in(0,\pi/2)$ such that, if a vector of length
$(1+\delta)^{-1}\mu^{-1}$ that forms the angle $\theta$ with a
given line is projected on that same line, then the length of this
projection is strictly greater than $(1+\alpha)\lambda \mu^{-1}$;
see Figure 1. Since
$$
(1+\delta)^{-1}\mu^{-1}-(1+\alpha)\lambda \mu^{-1}\geq
\lambda(1-\lambda)\mu^{-1}/4,
$$
we can for example pick $\theta$ such that

\begin{eqnarray*}
\cos\theta & = & \frac{(1+\alpha)\lambda
\mu^{-1}+\lambda(1-\lambda)\mu^{-1}/8}{(1+\delta)^{-1}\mu^{-1}}\\
& = & (1+\delta)[(1+\alpha)\lambda+\lambda(1-\lambda)/8].
\end{eqnarray*}

\begin{figure}
\centering \mbox{\subfigure[The inner
spiral.]{\epsfig{file=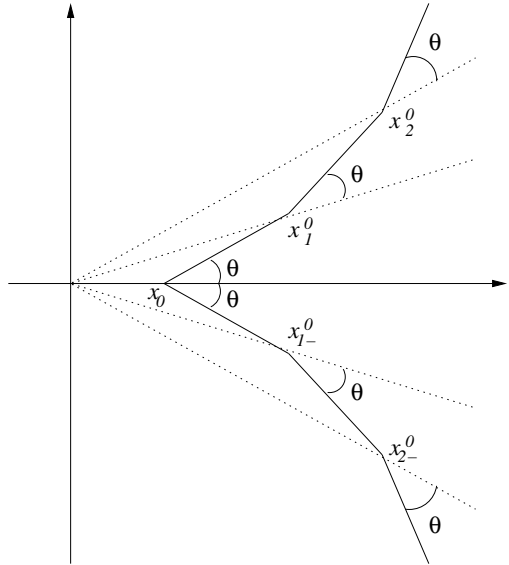,height=7.5cm}}}\par
\mbox{\subfigure[The branches hit $\partial B(0,u)$ at points
$\{y_i\}$. ]{\epsfig{file=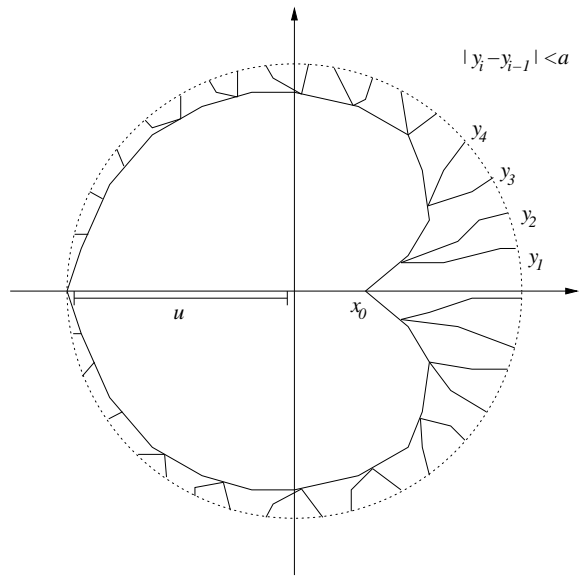,height=7.6cm}}}
\caption{Geometrical construction.}
\end{figure}

\noindent Now fix a point $x_0\in Q$ located on the positive
$x$-axis and draw two line segments starting from $x_0$ with angle
$\theta$ and $-\theta$ respectively to the $x$-axis; see Figure
2(a). The length of the segments is taken to be
$$
l=\frac{\varepsilon\lambda\mu^{-1}}{2(1+\alpha)}
$$
\noindent (the reason for this particular choice of $l$ will be
clear later on in the proof). Let $l_1^0$ and $l_{1-}^0$ denote
the two segments and write $x_1^0$ and $x_{1-}^0$ for the terminal
points (the zeroes in the superscripts will be explained later).
From the points $x_1^0$ and $x_{1-}^0$ we draw two new line
segments $l_2^0$ and $l_{2-}^0$ of length $l$. The segment $l_2^0$
($l_{2-}^0$) should form the angle $\theta$ ($-\theta$) with an
imaginary line through the origin and $x_1^0$ ($x_{1-}^0$). We
continue to draw line segments like this in an outward spiral
running in both directions. The segments should all be of length
$l$ and depending on the sign of the $y$-coordinate of its
starting point it should form the angle $\theta$ or $-\theta$ with
an imaginary line through the origin and its starting point.
Eventually the two spiral arms will meet at a point $(-u,0)$ on
the negative $x$-axis. Let $2n$ denote the number of segments
needed to achieve this. We then have two sets of line segments,
$\{l_k^0\}_1^n$ and $\{l_{k-}^0\}_1^n$, constituting the upper and
lower spiral arm respectively, and two sets of terminal points for
the line segments, $\{x_k^0\}_1^n$ and $\{x_{k-}^0\}_1^n$. Note
for the future that $|x_k^0|-|x_{k-1}^0|\geq l\cos\theta$, which
implies that

\begin{equation}\label{eq:stfor0}
|x_k^0|\geq (1+2\varepsilon)\lambda\mu^{-1}+kl\cos\theta.
\end{equation}

\noindent Now extend the construction by adding more edges, still
of length $l$, branching out from the points $\{x_k^0\}$ and
$\{x_{k-}^0\}$ towards the boundary of the circle with radius $u$
around the origin; see Figure 2(b). The branches should be built
up so that edges hit $\partial B(0,u)$ at points $\{y_j\}$ located
not more than a distance $a$ from each other, where
$$
a=\frac{1-\alpha}{1+\alpha}\cdot \frac{l}{4}
$$
\noindent (the choice of $a$ is motivated later). Furthermore, the
number of segments used to join a point $x_k^0$ (or $x_{k-}^0$) to
a point on the circle boundary $\partial B(0,u)$ should not exceed
$n-k$. We group the line segments in \emph{generations} depending
on how many links away from $x_0$ they are: An edge whose starting
point is linked to $x_0$ using $k-1$ other edges is placed in
generation $k$. If $l_k^i$ ($l_{k-}^i$) denotes segment number $i$
in generation $k$ in the upper (lower) half plane, we thus have
$n$ generations $\{l_k^i,l_{k-}^i\}_{i\geq 0}$, where the edges
with $i=0$ belongs to the inner spiral. Let $\{x_k^i\}$ and
$\{x_{k-}^i\}$ denote the terminal points of the segments $l_k^i$
and $l_{k-}^i$ respectively. The last demand on the construction
is that (\ref{eq:stfor0}) should hold for all terminal points in
generation $k$, that is,

\begin{equation}\label{eq:stfor}
\min_i|x_k^i|\geq (1+2\varepsilon)\lambda \mu^{-1}+kl\cos\theta.
\end{equation}

\noindent When the line segments are arranged, we complete the
construction by forming channels of width, say, $l/100$ around all
segments.\medskip

\noindent Write $m(x_0)$ for the total number of channels required
in the above construction. As indicated, this number depends on
the choice of the starting point $x_0$. Let $m$ denote the largest
value of $m(x_0)$ when $x_0\in Q$ and pick a time point $t_0$ that
fulfills the conditions (i)-(iii) described below. Some of these
conditions might seem awkward at first, but their purpose will
gradually become clear.

\begin{itemize}
\item[(i)] Write $\textbf{t}=(t,0,\ldots,0)$. Combining Lemma \ref{lemma:formbeg}
and Lemma \ref{lemma:mubtmu} yields that
$$
\lim_{b\rightarrow\infty}\lim_{t\rightarrow\infty}
\frac{\tilde{T}^b(\textbf{t})}{t}=\mu\quad \textrm{a.s.}
$$
Hence, for each $p>0$ and $\delta>0$, we have

\begin{equation}\label{eq:tubhastbeg}
P\left(\tilde{T}^b(\textbf{t})\geq (1+\delta)\mu t\right)\leq p
\end{equation}

\noindent if $b$ and $t$ are large. Let $t_0$ be large enough to
ensure that (\ref{eq:tubhastbeg}) holds for $b\geq t_0l/100$ and
$t\geq t_0l$ when $\delta=(1-\lambda)/2\lambda$ and $p=1/2m$.

\item[(ii)] Write $S_t$ for the infected area at time $t$ in a one-type
process with unit rate and let $t_0$ be large enough to guarantee
that
$$
P\left((1-\alpha)B\left(0,\mu^{-1}\right)\subset
\frac{S_t}{t}\subset
(1+\alpha)B\left(0,\mu^{-1}\right)\right)>1-\frac{1}{2m}
$$
for $t\geq t_0\varepsilon\lambda/8(1+\alpha)^2$.

\item[(iii)] Let $\xi=a/4u$ and $\delta=1/2$ in Lemma
\ref{lemma:andlutbrqu} and pick $t_0$ so that $t_0u\geq r_0^*$.
\end{itemize}

\noindent We will show that (\ref{eq:kritbeg}) holds for $t\geq
t_0$. To this end, fix $t\geq t_0$ and note that on $E_t$ we can
pick a point $x_0\in (\tilde{S}_t^2\cap tQ)/t$ to serve as
starting point for the geometrical construction described above
(by rotation invariance we may assume that $x_0$ is located on the
positive $x$-axis). At time $t$ we then have a strongly type 1
infected point $tx_0$ with $|x_0|\geq
(1+2\varepsilon)\lambda\mu^{-1}$. Also, since $\alpha\leq
\varepsilon$, on $F_t$ the type 2 infected area at time $t$ does
not reach further than $(1+\varepsilon)t\lambda \mu^{-1}$ from the
origin. Now define
$$
t'=(1+\delta)\mu tl
$$
\noindent and consider the state of the infection at time $t+t'$.
For the type 2 infection, by the choice of $l$ and $\alpha$ we
have on $F_t$ that

\begin{eqnarray}
\|S_{t+t'}^2\| & \leq & (1+\alpha)(t+t')\lambda
\mu^{-1}\nonumber\\
& \leq & \left(1+\frac{3}{2}\varepsilon\right)t\lambda\mu^{-1}
\label{eq:svagmatta1}.
\end{eqnarray}

\noindent To deal with the type 1 infection, let $\tilde{T}_{k}^i$
($\tilde{T}_{k-}^i$), $k\geq 1$, denote the time counting from
$t+(k-1)t'$ until the terminal point of the segment $l_k^i$
($l_{k-}^i$) is strongly type 1 infected assuming that at time
$t+(k-1)t'$ all type 1 infection is erased and replaced by a
$\gamma$-ball around the starting point of $l_k^i$ ($l_{k-}^i$)
while the type 2 infection is left as in the original process. By
(\ref{eq:svagmatta1}), on $F_t$ the type 2 infection has not yet
reached any parts of the channels in the first generation at time
$t+t'$. Hence, up to time $t+t'$ the spread of the type 1
infection inside the first generation channels behaves like
hampered one-type processes with $b=tl/100$. On the scale $t$ the
channels have length $tl$ and, since $t\geq t_0$, it follows from
the condition (i) in the choice of $t_0$ that on $E_t$ we have

\begin{equation}\label{eq:fgptbeg}
P^{1,\lambda}\left(\tilde{T}_1^i\geq
t'|\mathcal{F}_t,F_t\right)\leq \frac{1}{2m}\quad \textrm{for all
}i,
\end{equation}

\noindent where $\tilde{T}_1^i$ can also be replaced by
$\tilde{T}_{1-}^i$. For the state of the infection at time
$t+2t'$, a similar reasoning as for the time $t+t'$ yields that

\begin{equation}\label{eq:svagbel2}
\|S_{t+2t'}^2\|\leq (1+2\varepsilon)t\lambda\mu^{-1},
\end{equation}

\noindent implying that the analog of (\ref{eq:fgptbeg}) holds
also for the second generation passage times $\{\tilde{T}_2^i\}$
and $\{\tilde{T}_{2-}^i\}$. Now note that

\begin{eqnarray*}
(1+\alpha)t'\lambda \mu^{-1} & = &
(1+\delta)(1+\alpha)\lambda tl\\
& \leq & (1+\delta)[(1+\alpha)\lambda+\lambda(1-\lambda)/8]tl\\
& = & tl\cos\theta.
\end{eqnarray*}

\noindent Using the definitions of $t'$ and $l$ and the fact that
$\alpha\leq \varepsilon$, this implies that on $F_t$ we have

\begin{equation}\label{eq:svagbel}
\|S_{t+kt'}^2\|\leq
(1+2\varepsilon)t\lambda\mu^{-1}+(k-2)tl\cos\theta
\end{equation}

\noindent for $k\geq 2$ and thus, by (\ref{eq:stfor}), at time
$t+kt'$ the type 2 infection has not yet reached any parts of the
channels surrounding the line segments in the $k$th generation. Up
to time $t+kt'$, the spread of the type 1 infection inside the
$k$th generation channels hence behaves like hampered one type
processes with $b=tl/100$. It follows from the condition (i) in
the choice of $t_0$ that on $E_t$ the bound in (\ref{eq:fgptbeg})
holds also for $\{\tilde{T}_k^i\}$ and $\{\tilde{T}_{k-}^i\}$ with
$k\geq 3$ so that hence
$$
P^{1,\lambda}\left(\tilde{T}_k^i>t'|\mathcal{F}_t,F_t\right)\leq
\frac{1}{2m}\quad \textrm{for all }k\geq 1 \textrm{ and }i\geq 0,
$$
where $\tilde{T}_k^i$ can also be replaced by $\tilde{T}_{k-}^i$.
Let $C_t$ denote the event that no passage time in the system
exceed $t'$, that is,
$$
C_t=\bigcap_{k,i}\{\tilde{T}_k^i\leq t'\cap \tilde{T}_{k-}^i\leq
t'\}.
$$
\noindent Since there are at most $m$ channels in the system, on
$E_t$ we obtain

\begin{equation}\label{eq:pctbeg}
P^{1,\lambda}(C_t|\mathcal{F}_t,F_t)>\frac{1}{2}.
\end{equation}

\noindent Note that on $E_tF_tC_t$ all points $\{ty_j\}$ on the
boundary of the circle $B(0,tu)$ are strongly type 1 infected at
time $t+nt'$.\medskip

\noindent The next step is to use the one-type shape theorem to
show that with large probability the strong type 1 infection at
the points $\{ty_j\}$ will expand and create a connected type 1
layer around the type 2 infection. To this end, define
$$
t''=\frac{\varepsilon\lambda t}{8(1+\alpha)^2}
$$
and note that $(1+\alpha)t''\lambda\mu^{-1}\leq tl/4$.
Furthermore, it follows from (\ref{eq:stfor}) that
$(1+2\varepsilon)t\lambda\mu^{-1}+(n-1)tl\cos\theta\leq tu$.
Combined with some algebraic manipulation, this gives that on
$F_t$ we have

\begin{equation}\label{eq:svagvidslut}
\|S_{t+nt'+t''}^2\|\leq tu-\frac{tl}{4}.
\end{equation}

\noindent Now, for each $j$, assume that at time $t+nt'$ a new
process is started by reducing the type 1 infection to the
$\gamma$-ball around the point $ty_j$. More precisely, at time
$t+nt'$ all type 1 infection except the one in $B(ty_j,\gamma)$ is
erased while the type 2 infection is left unchanged. For $s\geq
t+nt'$, let $S_s^{1(j)}$ denote the type 1 infected region at time
$s$ in such a process and, for $s\geq 0$, define
$$
A^j_s=\left\{B\left(ty_j,(1-\alpha)s\mu^{-1}\right)\subset
S^{1(j)}_{t+nt'+s}\subset
B\left(ty_j,(1+\alpha)s\mu^{-1}\right)\right\}.
$$
Since $(1+\alpha)t''\mu^{-1}=tl/4$, the event $A^j_{t''}$ does not
depend on the state of the infection outside $B(ty_j,tl/4)$ and,
by (\ref{eq:svagvidslut}), $B(ty_j,tl/4)$ does not contain any
type 2 infection at time $t+nt'+t''$. Hence the one-type shape
theorem can be applied to estimate the probability of the event
$A^j_{t''}$ and since $t\geq t_0$ it follows from the condition
(ii) in the choice of $t_0$ that
$$
P^{1,\lambda}\left(A^j_{t''}\right)\geq 1-\frac{1}{2m}
$$
on $E_t$. Let
$$
A_t=\bigcap_jA^j_{t''}.
$$
The number of points $\{ty_j\}$ on $\partial B(0,tu)$ is clearly
bounded by $m$, implying that on $E_t$ we have

\begin{equation}\label{eq:pabeg}
P^{1,\lambda}(A_t|\mathcal{F}_t,F_t,C_t)>\frac{1}{2}.
\end{equation}

\noindent In words, $A_t$ is the event that all circles with
radius $(1-\alpha)t''\mu^{-1}$ around the points $\{ty_j\}$ are
type 1 infected at time $t+nt'+t''$. Since
$(1-\alpha)t''\mu^{-1}=ta$, where $ta$ is recognized as the
distance between the points $\{ty_j\}$, the circles overlap each
other so that a layer of type 1 infection with thickness at least
$ta/2$ concentrated around $\partial B(0,tu)$ is created.\medskip

\noindent Remember that the aim is to establish
(\ref{eq:kritbeg}). Trivially,
$$
P^{1,\lambda}(G_2^c|\mathcal{F}_t,F_t)\geq
P^{1,\lambda}(G_2^c|\mathcal{F}_t,F_t,C_t,A_t)\cdot
P^{1,\lambda}(A_t|\mathcal{F}_t,F_t,C_t)\cdot
P^{1,\lambda}(C_t|\mathcal{F}_t,F_t)
$$
and using (\ref{eq:pctbeg}) and (\ref{eq:pabeg}) it follows that
on $E_t$ we have

\begin{equation}\label{eq:tupp}
P^{1,\lambda}(G_2^c|\mathcal{F}_t,F_t)>P^{1,\lambda}
(G_2^c|\mathcal{F}_t,F_t,C_t,A_t) \cdot 2^{-2}.
\end{equation}

\noindent What remains is to bound the probability that $G_2^c$
occurs on $E_tF_tC_tA_t$ from below. To do this, note that on
$E_tF_tC_tA_t$, at time $t+nt'+t''$ the type 2 infection is
contained in $B(0,tu)$ and the annulus $B(0,tu+ta/4)\backslash
B(0,tu)$ is filled with type 1 infection. In between the type 1
layer and the type 2 infection there might however still be
uninfected regions. Clearly we are done if we can find a lower
bound for the probability that $G_2^c$ occurs when these regions
are assumed to be occupied by type 2 infection. Hence, consider a
two-type growth process with infection rates $(1,\lambda)$ started
from a connected configuration without holes such that
$S_0^2\subset B(0,tu)$ and $S_t^1\cup S_t^2\supset B(0,tu+ta/4)$.
It follows from Lemma \ref{lemma:andlutbrqu} and the condition
(iii) in the choice of $t$ that the probability that the type 2
infection never makes an effective outburst in such a process is
at least 1/2. Combining this with (\ref{eq:tupp}) yields
$$
P^{1,\lambda}(G_2^c|\mathcal{F}_t,F_t)>\frac{1}{2^3}.
$$
Hence (\ref{eq:kritbeg}) is established for $d=2$ and the
proposition is proved.\medskip

\noindent For $d\geq 3$, the geometrical construction is obtained
by first rotating the two-dimensional inner spiral around the
$x$-axis a finite number of times and then add branches --
emanating from the rotated spiral arms -- that hit the surface of
the ball $B(0,u)$ closely enough. As in the two-dimensional case,
Lemma \ref{lemma:formbeg} and Lemma \ref{lemma:mubtmu} can be
combined to show that, if $t$ is large, then with large
probability the type 1 infection travels fast enough through the
channels to reach the points on $\partial B(0,u)$ in time to
create a thick layer around the type 2 infection. The rest of the
proof is analogous.$\hfill\Box$\medskip

\section{Proof of Theorem \ref{th:huvudres}}

\noindent In this section we prove Theorem \ref{th:huvudres} using
arguments similar to those used for the main result of
H\"aggstr\"om and Pemantle (2000), but with a different twist at
the end, which is needed because of the unboundedness of the
outbursts radii. The proof is based on Proposition \ref{prop:key}
and a coupling of the two-type processes with distributions
$\{P^{1,\lambda}\}_{\lambda\geq 0}$ valid for all
$\lambda\in[0,1]$ simultaneously.\medskip

\noindent The part of the proof where Proposition \ref{prop:key}
comes into play is formulated separately in the following lemma,
which says roughly that the shape theorem holds on the event of
unbounded growth for the weaker infection type and that the radius
of the asymptotic shape in this case is determined by the weaker
infection type.

\begin{lemma}\label{lemma:keylemma} Let $S_t^1\cup S_t^2$ be the
region infected at time $t$ in a two-type process with
distribution $P^{1,\lambda}$, where $\lambda\in [0,1]$, and assume
that (\ref{eq:rvillk}) holds. Then, for any
$\varepsilon\in(0,\lambda\mu^{-1})$, we have $P^{1,\lambda}$-a.s.
on the event $G_2$ that
$$
(1-\varepsilon)B\left(0,\lambda\mu^{-1}\right)\subset
\frac{S_t^1\cup S_t^2}{t}\subset
(1+\varepsilon)B\left(0,\lambda\mu^{-1}\right)
$$
for all sufficiently large $t$.
\end{lemma}

\noindent \emph{Proof:} Let $\|S_t^1\cup S_t^2\|_*$ denote the
minimum distance from the origin to $(S_t^1\cup S_t^2)^c$, that
is,
$$
\|S_t^1\cup S_t^2\|_*=\sup\{s;\hspace{0.1cm} B(0,s)\subset
S_t^1\cup S_t^2\}.
$$
The lemma follows if we can show that
$$
\frac{\|S_t^1\cup S_t^2\|}{t}\rightarrow \lambda\mu^{-1}\quad
\textrm{and}\quad \frac{\|S_t^1\cup S_t^2\|_*}{t}\rightarrow
\lambda\mu^{-1}
$$
$P^{1,\lambda}$-a.s. on the event $G_2$. Since $\|S_t^1\cup
S_t^2\|_*\leq \|S_t^1\cup S_t^2\|$ it suffices to prove that
$P^{1,\lambda}$-a.s. on $G_2$ we have

\begin{equation}\label{eq:kl1}
\limsup_{t\rightarrow \infty}\frac{\|S_t^1\cup S_t^2\|}{t}\leq
\lambda \mu^{-1} \end{equation}

\noindent and

\begin{equation}\label{eq:kl2}
\liminf_{t\rightarrow \infty}\frac{\|S_t^1\cup S_t^2\|_*}{t}\geq
\lambda \mu^{-1}.
\end{equation}

\noindent The lower bound (\ref{eq:kl2}) follows immediately from
Lemma \ref{lemma:ttpetp} and Theorem \ref{th:aformgen}. To
establish (\ref{eq:kl1}), note that, by Lemma \ref{lemma:svagdom}
and Theorem \ref{th:aformgen},
$$
\limsup_{t\rightarrow \infty}\frac{\|S_t^2\|}{t}\leq \lambda
\mu^{-1}.
$$
We are done if we can show that $S_t^2$ can also be replaced by
$S_t^1$ here, that is, if we can show that

\begin{equation}\label{eq:kl3}
\limsup_{t\rightarrow \infty}\frac{\|S_t^1\|}{t}\leq \lambda
\mu^{-1}.
\end{equation}

\noindent But this is a consequence of Proposition \ref{prop:key},
since, if (\ref{eq:kl3}) fails, there is an $\varepsilon>0$ such
that the type 1 infected region reaches outside
$(1+3\varepsilon)tB(0,\lambda \mu^{-1})$ for arbitrarily large $t$
and by Proposition \ref{prop:key} this prevents the event
$G_2$.$\hfill \Box$\medskip

\noindent Moving on to the aforementioned simultaneous coupling of
the two-type processes with distributions
$\{P^{1,\lambda}\}_{\lambda\in[0,1]}$, let $N_1$ and $N_2$ be two
independent unit rate Poisson processes. We will couple the growth
processes by successively thinning the Poisson process $N_2$ and
then use it to generate the type 2 outbursts. This is done as
follows: Associate independently to each point in $N_2$ a random
variable uniformly distributed over $[0,1]$, and let $\lambda N_2$
be the set of points in $N_2$ whose attached uniform variable is
smaller than or equal to $\lambda$. Then $\lambda N_2$ is a
Poisson process with rate $\lambda$ and hence, for each
$\lambda\in[0,1]$, a two-type process $\{S_t^1(\lambda)\cup
S_t^2(\lambda)\}_{t\geq 0}$ with distribution $P^{1,\lambda}$ is
obtained by starting from $B(-\textbf{2}\gamma,\gamma)$ and
$B(0,\gamma)$ at time zero and then using $N_1$ to generate the
type 1 outbursts and $\lambda N_2$ to generate the type 2
outbursts. Write $Q$ for the probability measure underlying this
coupling and let $G_i^{\lambda}$ denote the event that the type
$i$ infection grows indefinitely at parameter value
$\lambda$.\medskip

\noindent \emph{Proof of Theorem \ref{th:huvudres}:} By
time-scaling and symmetry we have
$$
P^{1,\lambda}(G)=P^{1,1/\lambda}(G)
$$
and hence it is enough to prove that $P^{1,\lambda}(G)=0$ for all
but at most countably many $\lambda\in[0,1]$. To this end, we will
show that, for any $\lambda'<\lambda\in [0,1]$, we have

\begin{equation}\label{eq:bh1}
Q(G_1^{\lambda}\cap G_2^{\lambda'})=0.
\end{equation}

\noindent This implies that with $Q$-probability 1 the event
$G_1^\lambda\cap G_2^\lambda$ occurs for at most one
$\lambda\in[0,1]$: By construction of the probability measure $Q$,
the event $G_1^\lambda$ is decreasing in $\lambda$ -- that is, if
$G_1^{\lambda}$ occurs then $G_1^{\lambda'}$ occurs for all
$\lambda'<\lambda$ -- and the event $G_2^\lambda$ is increasing in
$\lambda$. Hence the set of lambdas for which the event
$G_1^\lambda\cap G_2^\lambda$ occurs is $Q$-a.s.\ an interval. If
with positive $Q$-probability the interval were non-degenerated,
there would be $\lambda'<\lambda$ in $[0,1]$ such that the event
$G_1^\lambda\cap G_2^\lambda\cap G_1^{\lambda'}\cap
G_2^{\lambda'}$ has positive $Q$-probability. This however
contradicts (\ref{eq:bh1}). Thus with $Q$-probability 1 the
interval consists of a single point, implying that $Q$-a.s.\ the
event $G_1^\lambda\cap G_2^\lambda$ occurs for at most one
$\lambda\in[0,1]$. Clearly
$$
P^{1,\lambda}(G)=Q(G_1^\lambda\cap G_2^\lambda)
$$
and hence it follows that $\{\lambda\in[0,1];\hspace{0.1cm}
P^{1,\lambda}(G)>0\}$ is countable.\medskip

\noindent To establish (\ref{eq:bh1}), fix
$\lambda'<\lambda\in[0,1]$ and assume that $G_2^{\lambda'}$
occurs. By Lemma \ref{lemma:keylemma} we then have
$$
\limsup_{t\rightarrow\infty}\frac{\|S_t^1(\lambda')\cup
S_t^2(\lambda')\|}{t}\leq \lambda'\mu^{-1}
$$
so that, in particular,
$$
\limsup_{t\rightarrow\infty}\frac{\|S_t^1(\lambda')\|}{t}\leq
\lambda'\mu^{-1}.
$$
From the construction of the $Q$-coupling it is clear that
$\|S_t^1(\lambda)\|\leq \|S_t^1(\lambda')\|$ and hence it follows
that

\begin{equation}\label{eq:bh2}
\limsup_{t\rightarrow\infty}\frac{\|S_t^1(\lambda)\|}{t}\leq
\lambda'\mu^{-1}.
\end{equation}

\noindent Furthermore, by Lemma \ref{lemma:keylemma}, if $t$ is
large

\begin{equation}\label{eq:bh3}
\limsup_{t\rightarrow\infty}\frac{\|S_t^1(\lambda)\cup
S_t^2(\lambda)\|_*}{t}\geq \lambda\mu^{-1}.
\end{equation}

\noindent Now pick $\varepsilon>0$ such that
$(1+\varepsilon)\lambda'<(1-\varepsilon)\lambda$. Combining
(\ref{eq:bh2}) and (\ref{eq:bh3}) yields that there is a time $T$
such that, for $t\geq T$, we have

\begin{equation}\label{eq:torsdag}
S_t^1(\lambda)\subset (1+\varepsilon)tB(0,\lambda'\mu^{-1})
\end{equation}

\noindent and

\begin{equation}\label{eq:fredag}
(1-\varepsilon)tB(0,\lambda\mu^{-1})\subset S_t^1(\lambda)\cup
S_t^2(\lambda).
\end{equation}

\noindent By (\ref{eq:fredag}), an outburst that occurs at a time
point $t\geq T$ must reach outside
$(1-\varepsilon)tB(0,\lambda\mu)$ to be effective and, by the
choice of $\varepsilon$, we have
$$
(1+\varepsilon)B(0,\lambda'\mu^{-1})\subset
(1-\varepsilon)B(0,\lambda\mu^{-1}).
$$
Thus an effective type 1 outburst at a time point $t\geq T$ would
cause the type 1 infected region to reach outside
$(1+\varepsilon)B(0,\lambda'\mu^{-1})$. This conflicts with
(\ref{eq:torsdag}) and hence no effective type 1 outbursts can
occur after time $T$. Clearly this prevents the event
$G_1^\lambda$. $\hfill\Box$

\section{Proof of Proposition \ref{prop:startomr}}

\noindent This section is devoted to the proof of Proposition
\ref{prop:startomr}.\medskip

\noindent\emph{Proof of Proposition \ref{prop:startomr}:} Pick
bounded sets $\Gamma_1,\Gamma_2,\Gamma_1',\Gamma_2'$ of positive
Lebesgue measure such that $\Gamma_1$ and $\Gamma_2$ and also
$\Gamma_1'$ and $\Gamma_2'$ are disjoint. We will show that if $G$
has positive probability in the process started from
$(\Gamma_1,\Gamma_2)$, then $G$ occurs with positive probability
in the process started from $(\Gamma_1',\Gamma_2')$ as well. To
this end, first consider the process started from
$(\Gamma_1,\Gamma_2)$. By Lemma \ref{lemma:andlutbr}(a) almost
surely only finitely many effective outbursts occur in the set
$\Gamma_1'\cup \Gamma_2'$ during the progress of the growth in
this process and hence there is a time $t<\infty$ such that with
probability, say, 1/2 no effective outbursts occur in
$\Gamma_1'\cup \Gamma_2'$ after time $t$. Let $U_i$ denote the set
of effective type $i$ outbursts that occur in the set $S_t^1\cup
S_t^2$ after time $t$. A second application of Lemma
\ref{lemma:andlutbr}(a) yields that the sets $U_i$ are almost
surely finite.\medskip

\noindent Now consider a process started from
$(\Gamma_1',\Gamma_2')$, coupled with the one started from
$(\Gamma_1,\Gamma_2)$ in such a way that the same Poisson
processes are used to generate the outbursts after time $t$.
Before time $t$ the process evolves independently of the one
started from $(\Gamma_1,\Gamma_2)$. We will describe a scenario
for this process that causes the infection to develop in the same
way as in the process started from $(\Gamma_1,\Gamma_2)$ after
time $t$. To prepare for this, join each point in $U_i$ by a curve
with a point in the interior of $\Gamma_i'$. The connections are
made by aid of concatenations of straight line segments and the
restrictions on a connection joining a point in $U_1$ ($U_2$) with
a point in $\Gamma_1'$ ($\Gamma_2'$) are:

\begin{itemize}
\item[--] It is not allowed to cross any part of $\Gamma_2'$ ($\Gamma_1'$).
\item[--] It should stay within the region infected at time $t$
in the process started from $(\Gamma_1,\Gamma_2)$.
\item[--] It can not pass through points in $U_2$ ($U_1$).
\end{itemize}

\noindent The first restriction might not be possible to fulfill
if the set $\Gamma_1'$ ($\Gamma_2'$) is enclosed by $\Gamma_2'$
($\Gamma_1'$). However, if this should be the case, we condition
on a large outburst occurring in $\Gamma_1'$ ($\Gamma_2'$) at some
early time point transmitting the type 1 (2) infection past
$\Gamma_2'$ ($\Gamma_1'$). Then we use the outer parts of the type
1 (2) infected region as terminal for the connections. When the
connections are made we let each one of them be surrounded by a
path of width $2\varepsilon$, where $\varepsilon>0$ is chosen
small enough to guarantee that the paths are disjoint. (In two
dimensions it is sometimes impossible to avoid paths from crossing
each other and hence we have to allow overlapping paths at
crossing points.) Let $\{P_i^k;\hspace{0.1cm} i=1,2 \textrm{ and
}k\geq 1\}$ denote the paths and write $T(P_i^k)$ for the time it
takes for the infection to wander along $P_i^k$ from $\Gamma_i'$
to its terminal point in $U_i$ by aid of $\varepsilon$-small
outbursts not reaching outside the path (in the two-dimensional
case we allow for outbursts with radius $2\varepsilon$ at the
possible crossings). Furthermore, define $\tau$ to be the time
when all points in $U_1\cup U_2$ are reached by the infection
using the paths, that is,
$$
\tau=\max_{i,k}\{T(P_i^k)\}.
$$
The desired scenario for the process started from
$(\Gamma_1',\Gamma_2')$ is now obtained as follows:

\begin{itemize}
\item[1.] Assume that the infection wanders along the paths from
$\Gamma_1'$ and $\Gamma_2'$ to the points in $U_1$ and $U_2$ by
aid of $\varepsilon$-small outbursts.

\item[2.] Suppose that $\tau\leq t$. Also assume that no outbursts except
for the ones on the paths occur before time $\tau$ and that no
outbursts at all occur in the time interval $(\tau,t)$. At time
$t$ then, the infected region consists of the initial sets
$\Gamma_1'$ and $\Gamma_2'$ together with fine infected strings
linking these sets to the points in $U_1$ and $U_2$.

\item[3.] After time $t$ the same Poisson processes as in the
process started from $(\Gamma_1,\Gamma_2)$ are used to generate
the outbursts. Hence we know that effective outbursts of the same
type as in the process started from $(\Gamma_1,\Gamma_2)$ will
take place at the points in $U_1$ and $U_2$. During the progress
of the growth it might happen that some parts of the region that
is infected at time $t$ in the process started from
$(\Gamma_1,\Gamma_2)$ are infected by another infection type.
Assume that no effective outbursts take place in those regions. By
Lemma \ref{lemma:andlutbr}(b) this event has positive probability.
\end{itemize}

\noindent Write $G_{\Gamma_1,\Gamma_2}$ for the event that both
infection types grow indefinitely in the process started from
$(\Gamma_1,\Gamma_2)$ and write $\hat{G}_{\Gamma_1',\Gamma_2'}$
for the same event in the coupled process started from
$(\Gamma_1',\Gamma_2')$. Trivially
$$
P(\hat{G}_{\Gamma_1',\Gamma_2'})\geq
P(\hat{G}_{\Gamma_1',\Gamma_2'}|G_{\Gamma_1,\Gamma_2})
P(G_{\Gamma_1,\Gamma_2}) .
$$
The second factor on the right-hand side is positive by
assumption. As for the first factor, note that if both infection
types grow indefinitely in the process started from
$(\Gamma_1,\Gamma_2)$, the above scenario guarantees mutual
unbounded growth also in the coupled process started from
$(\Gamma_1',\Gamma_2')$, since in both processes the only
outbursts that will reach outside the region infected at time $t$
in the process started from $(\Gamma_1,\Gamma_2)$ are the ones in
$U_1$ and $U_2$. Furthermore, the above scenario clearly has
positive probability because it only depends on finitely many
outbursts. Hence also the first factor is positive and it follows
that $P(\hat{G}_{\Gamma_1',\Gamma_2'})>0$.  Since
$P(\hat{G}_{\Gamma_1',\Gamma_2'})=
P_{\Gamma_1',\Gamma_2'}^{\lambda_1,\lambda_2}(G)$, we are done.
\hfill $\Box$

\section*{References}

\noindent Biggins, J.D. (1995): The growth and spread of the
general branching random walk, \emph{Ann. Appl. Probab.}, vol.
\textbf{5}, 1008-1024.\medskip

\noindent Deijfen, M. (2003): Asymptotic shape in a continuum
growth model, \emph{Adv. Appl. Probab.}, vol.\textbf{35:2},
303-318 .\medskip

\noindent H\"{a}ggstr\"{o}m, O. and Pemantle, R. (1998): First
passage percolation and a model for competing spatial growth,
\emph{J. Appl. Probab.}, vol. \textbf{35}, 683-692.\medskip

\noindent H\"{a}ggstr\"{o}m, O. and Pemantle, R. (2000): Absence
of mutual unbounded growth for almost all parameter values in the
two-type Richardson model, \emph{Stoch. Proc. Appl.}, vol.
\textbf{90}, 207-222.\medskip

\noindent Jagers, P. (1975): \emph{Branching Processes with
Biological Applications}, Wiley.\medskip

\noindent Richardson, D. (1973): Random growth in a tessellation,
\emph{Proc. Cambridge Phil. Soc.} \textbf{74}, 515-528.\medskip

\noindent Williams, D. (1991): \emph{Probability with
Martingales}, Cambridge University Press.\medskip

\pagebreak

\begin{center}\textbf{Corrections}\end{center}

\noindent $\bullet$ At the end of the proof of Proposition 2.1 in
Paper II, page 225, it should be noted that, by symmetry, the same
convergence result as for the position of the rightmost individual
holds for the position of the leftmost one as well. ``Rightmost''
(``leftmost'') refers to the maximal (minimal) element in the set
of positions.\medskip

\noindent $\bullet$ The notation $\|\cdot\|$ is defined on page
235, line 22, but appears for the first time on page 230.\medskip

\noindent $\bullet$ The choice of the angle $\theta$ on page 237
is inappropriate. This is adjusted as follows.

\begin{itemize}
\item[--] The right-hand side on line 7, page 237, should be
$\lambda(1-\lambda)\mu^{-1}/4$.
\item[--] The denominator in the fraction on line 9, page 237, is set to $(1+\alpha)\lambda\mu^{-1}+
\lambda(1-\lambda)\mu^{-1}/8$ and the final expression on the same
line hence becomes
$(1+\delta)[(1+\alpha)\lambda+\lambda(1-\lambda)/8]$.
\item[--] Page 241, line 1:
``$\leq(1+\delta)[(1+\alpha)\lambda+\lambda(1-\lambda)/8]tl=tl\cos\theta$''.
\end{itemize}

\noindent $\bullet$ The inequalities (5.10) and (5.13) are
incorrectly motivated. This is adjusted as follows.

\begin{itemize}
\item[--] On page 241, the lines 2 to 4 and -10 (from ``Furthermore..'') to -5 are
deleted.
\item[--] Page 241, line 2:  ``Using the definitions of $t'$ and $l$ and the fact that
$\alpha\leq \varepsilon$, this implies that on $F_t$ we have''.
\item[--] Page 242, lines 1 and 2: ``and note that $(1+\alpha)t''\lambda\mu^{-1}
\leq tl/4$. Furthermore, it follows from (16) that
$(1+2\varepsilon)t\lambda\mu^{-1}+(n-1)tl\cos\theta\leq tu$.
Combined with some algebraic manipulation, this gives that on
$F_t$ we have''.
\end{itemize}

\noindent Finally we list a few misprints:

\begin{tabbing}
\textbf{page}\= \hspace{0.6cm}\textbf{line}\hspace{0.6cm}\=
\textbf{reads}\hspace{3.5cm}\=
\textbf{should read}\\
\hspace*{0.3cm}224\> \hspace{0.95cm}-4\> $R$\> $2R$\\
\hspace*{0.3cm}225\> \hspace{0.95cm}-8\> $\zeta$\> $2\zeta$\\
\hspace*{0.3cm}231\> \hspace{1.05cm}9\> $\leq$\> $=$\\
\hspace*{0.3cm}233\> \hspace{0.9cm}19\> $S_0$\> $S^*_0$\\
\hspace*{0.3cm}234\> \hspace{0.9cm}11\> $r^{d+1}$\> $r$\\
\hspace*{0.3cm}235\> \hspace{0.8cm}-10\> $B(0,\lambda\mu^{-1})$\> $\lambda\mu^{-1}$\\
\hspace*{0.3cm}236\> \hspace{0.9cm}17\> $\|S_s^2\|$\> $S_s^2$\\
\hspace*{0.3cm}237\> \hspace{0.95cm}-5\> $x$-coordinate\> $y$-coordinate\\
\hspace*{0.3cm}239\> \hspace{0.5cm}-5,-3\> $T^b(\textbf{t})$\> $\tilde{T}^b(\textbf{t})$\\
\hspace*{0.3cm}242\> \hspace{1.05cm}3\> $u$\> $tu$\\
\hspace*{0.3cm}248\> \hspace{0.9cm}11\> (2002)\> (2003)\\
\end{tabbing}

\end{document}